%% file: hmfa-openai.tex
\documentclass[11pt,onecolumn]{article}
\usepackage[top=1 in, bottom=1 in, left=0.75 in, right=0.75 in]{geometry}
\usepackage{hyperref}
\usepackage{amsmath} 
\usepackage{amssymb}
\usepackage{tabularx}
\usepackage[sc,osf,noBBpl]{mathpazo} 
\usepackage{cite} 
\usepackage{graphicx} 
\usepackage{subcaption}
\usepackage{algorithm}
\usepackage{algorithmic}
\usepackage{amssymb}
\usepackage{epsfig,xcolor,url}
\usepackage{epstopdf}
\usepackage{bbm}
\usepackage{enumitem} 

\usepackage{stfloats}
\usepackage{tikz}
\usetikzlibrary{arrows.meta,positioning}

\usepackage[utf8]{inputenc}
\usepackage{amsthm,bm} 
\renewenvironment{proof}{{\sffamily\bfseries Proof. }}{\qed}
\newtheorem{theorem}{Theorem}
\newtheorem{lemma}{Lemma}
\newtheorem{proposition}{Proposition}
\newtheorem{corollary}{Corollary}

\theoremstyle{definition}

\newtheorem{assumption}{Assumption}
\newtheorem{remark}{Remark}

\usepackage{thmtools}
\usepackage{thm-restate} 
\usepackage{enumitem}

\title{Second-Order Steady-State Approximation of Heterogeneous Mean-Field Models}

\author{
	Lei Ying \\
	University of Michigan, Ann Arbor \\  \texttt{leiying@umich.edu}\\
}

\begin{document}
\date{}
\maketitle

\begin{abstract}
Mean-field models provide tractable approximations for large stochastic systems, but second-order steady-state results are largely limited to homogeneous or fixed-dimensional population processes. This paper develops a second-order steady-state approximation for fully heterogeneous mean-field models. Conditioning on the realized agent parameters, we study fixed-dimensional weighted occupancies even though the dimension of the full heterogeneous mean-field system grows linearly with the number of agents $M$. We first establish an $O(M^{-1})$ bound on the mean-field approximation error for stationary expectations of twice continuously differentiable functions of these occupancies. We then characterize the leading term in this error in terms of the stationary mean bias and covariance of the weighted occupancies. A Lyapunov equation for the linearized dynamics yields an approximation to the stationary covariance of the unscaled weighted occupancy with error $O(M^{-3/2})$. Its solution, together with the drift derivatives at equilibrium, also gives an explicit $1/M$ correction to the mean-field prediction of the stationary mean. The corrected mean has error $o(M^{-1})$, which improves to $O(M^{-3/2})$ under a stronger bound on the drift remainder. The resulting covariance and mean-refinement equations have the same structure as the finite-time heterogeneous refined mean-field equations, but our results directly provide steady-state guarantees for weighted occupancies rather than marginal state probabilities of individual agents. The results are illustrated through a heterogeneous power-of-two load-balancing model with buffer size two.
\end{abstract}

\section*{\normalsize The Role of Generative AI (GenAI) }
This paper is a collaborative effort between the author and GenAI. ChatGPT was extensively used in writing this manuscript. Most proofs were generated, corrected, and refined through iterative
interactions between ChatGPT and the author. The author verified all results, organized
and edited the paper, and takes full responsibility for the
final content of the manuscript.

\input{intro}

\input{models}

\input{wfo}
\input{so}

\input{po2}

\section{Conclusion}

We characterize the steady-state covariance and bias of weighted occupancies
in fully heterogeneous mean-field models conditional on the realized agent
parameters. The matrix $\boldsymbol{\Gamma}_\varphi/M$ obtained from the
Lyapunov equation approximates the occupancy covariance with error
$O(M^{-3/2})$. Under the additional drift condition A9, we obtain an
explicit $1/M$ correction for stationary expectations of twice continuously
differentiable functions of the occupancy, with residual error $o(M^{-1})$.
We verify the assumptions for heterogeneous power-of-two load balancing
with buffer size two.

\bibliographystyle{plain}
\bibliography{inlab-refs}
\end{document}

%% file: intro.tex
\section{Introduction}
\label{sec:introduction}

Mean-field models provide tractable deterministic approximations for large stochastic systems and have been widely used in epidemic models, queueing and load-balancing systems, communication networks, and interacting particle systems. Consider $M$ agents evolving on a common $n$-state space. Classical finite-time convergence results go back to Kurtz~\cite{Kur_71,Kur_81}. Using Stein's method, Ying~\cite{Yin_16} established an $O(1/M)$ bound on the steady-state mean-square approximation error for a broad class of homogeneous mean-field models. More recently, Ying~\cite{Yin_26} extended this result to fully heterogeneous systems under assumptions that include stability and uniform robustness of the occupancy dynamics.

This paper studies the second-order steady-state behavior of fully heterogeneous mean-field models. For homogeneous density-dependent population processes, Gast and Van Houdt~\cite{GasHou_17} developed a refined mean-field approximation in steady state, in which the leading $1/M$ correction is characterized by the Jacobian and Hessian of the mean-field drift and the solution of a Lyapunov equation. Barbour, Luczak, and Xia~\cite{BarLucXia_18a,BarLucXia_18b} developed quantitative multivariate discrete-normal approximations for equilibrium distributions of density-dependent Markov jump processes. These results rely on a fixed-dimensional population-process representation. For fully heterogeneous systems, Allmeier and Gast~\cite{AllGas_22} established mean-field and refined mean-field approximation results over finite time horizons. These finite-time results leave open the corresponding steady-state question when the dimension of the full mean-field system grows with the number of agents.

We address this problem by conditioning on the realized agent parameters and studying weighted occupancies. Let $\theta_m$ denote the parameter of agent $m$, $\mathbf{W}_m(t)\in\mathbb{R}^{n-1}$ its reduced one-hot state, and $\mathbf{w}_m^*$ its component of the quenched mean-field equilibrium. For a fixed bounded weight function $\phi$, define $\mathbf{X}_\phi(t)=M^{-1}\sum_{m=1}^M\phi(\theta_m)\mathbf{W}_m(t)$, $\mathbf{x}_\phi^*=M^{-1}\sum_{m=1}^M\phi(\theta_m)\mathbf{w}_m^*$, and $\mathbf{Y}_\phi=\sqrt{M}\bigl(\mathbf{X}_\phi(\infty)-\mathbf{x}_\phi^*\bigr)$, where $\mathbf{X}_\phi(\infty)$ denotes the stationary weighted occupancy. We linearize the heterogeneous mean-field system at its equilibrium and let $\overline{\boldsymbol{\Sigma}}$ denote the solution of the corresponding Lyapunov equation. With $\boldsymbol{\Gamma}_\phi$ denoting its weighted projection, we show, under assumptions A1--A8, that
\begin{align}
\bigl\|\operatorname{Cov}(\mathbf{Y}_\phi\mid\boldsymbol{\theta})-\boldsymbol{\Gamma}_\phi\bigr\|_{\mathrm{op}}
&= O(M^{-1/2}). \label{eq:intro-weighted-covariance}
\end{align}
Here $\boldsymbol{\theta}=(\theta_1,\ldots,\theta_M)$ is the realized parameter vector. Thus, although the full mean-field dimension grows linearly with $M$, the covariance matrix in~\eqref{eq:intro-weighted-covariance} has size $(n-1)\times(n-1)$, independent of $M$. Equivalently, $\boldsymbol{\Gamma}_\phi/M$ approximates the covariance of the unscaled occupancy $\mathbf{X}_\phi(\infty)$ with error $O(M^{-3/2})$.

We further characterize the leading finite-size bias. Under the additional drift condition A9, any twice continuously differentiable function $h:\mathbb{R}^{n-1}\to\mathbb{R}$ satisfies
\begin{align}
\mathbb{E}\bigl[h(\mathbf{X}_\phi(\infty))\mid\boldsymbol{\theta}\bigr]
&=h(\mathbf{x}_\phi^*)+\frac{1}{M}\left(Dh(\mathbf{x}_\phi^*)[\boldsymbol{\beta}_\phi]
+\frac{1}{2}\operatorname{Tr}\!\left(D^2h(\mathbf{x}_\phi^*)\boldsymbol{\Gamma}_\phi\right)\right)+o(M^{-1}). \label{eq:intro-weighted-performance}
\end{align}
Here $Dh$ and $D^2h$ denote the first and second derivatives of $h$, respectively, $\operatorname{Tr}$ denotes the matrix trace, and $\boldsymbol{\beta}_\phi\in\mathbb{R}^{n-1}$ is the mean-refinement vector for the weighted occupancy, determined by the Jacobian and Hessian of the mean-field drift and $\overline{\boldsymbol{\Sigma}}$. The refined stationary mean $\mathbf{x}_\phi^*+\boldsymbol{\beta}_\phi/M$ has error $o(M^{-1})$, which improves to $O(M^{-3/2})$ under the stronger drift remainder bound in A9. In the homogeneous density-dependent setting with ordinary occupancy
($\phi\equiv1$), the covariance approximation in~\eqref{eq:intro-weighted-covariance}
and the $1/M$ performance correction in~\eqref{eq:intro-weighted-performance}
have the same form as in the classical steady-state refinement in~\cite{GasHou_17}. The distributional approximations in~\cite{BarLucXia_18a,BarLucXia_18b} also provide bounds in total variation. Here we obtain the covariance and $1/M$ performance refinement directly in the fully heterogeneous, growing-dimensional setting.

The covariance and mean-refinement equations have the same structure as the finite-time refined mean-field equations in~\cite{AllGas_22}. In~\cite{AllGas_22}, the covariance correction satisfies a time-dependent Lyapunov equation, while the mean correction is driven by the Hessian of the mean-field drift contracted with this covariance. Formally setting the time derivatives to zero at the mean-field equilibrium yields the equations defining $\overline{\boldsymbol{\Sigma}}$ and the mean-refinement vector in this paper. The main theoretical results in~\cite{AllGas_22} concern finite-time marginal state probabilities of individual heterogeneous agents. As discussed in~\cite{Yin_26}, dependence among agents can prevent the occupancy measure from concentrating even when the individual marginal probabilities are approximated accurately. We directly characterize the steady-state covariance and finite-size bias of weighted occupancies.

As a preliminary step toward these second-order results, we extend the first-order steady-state result of~\cite{Yin_26} from the ordinary occupancy to weighted occupancies. This provides the second-moment and bias estimates needed in the second-order analysis. We use a local coupling argument to establish a fourth-moment bound and give structural conditions on the drift and instantaneous covariance rate that can be checked from the transition rules and their smooth extensions. We verify the assumptions for heterogeneous power-of-two load balancing with buffer size two.

%% file: models.tex
\section{Models}
\label{sec:models}

\subsection{Heterogeneous $M$-Agent CTMC}
\label{subsec:heterogeneous-ctmc}

Consider a continuous-time Markov chain (CTMC) consisting of $M$
heterogeneous agents. Each agent $m$ is associated with a parameter
$\theta_m$ belonging to a bounded set $\Theta$. Throughout the paper, we
analyze the system conditional on the realized parameter vector
$\boldsymbol\theta=(\theta_1,\ldots,\theta_M)$.

Each agent evolves on the common finite state space
${\cal S}=\{0,\ldots,n-1\}$, where $n$ is independent of $M$. Let
$\boldsymbol e_1,\ldots,\boldsymbol e_{n-1}$ denote the standard basis
vectors of $\mathbb R^{n-1}$ and let $\boldsymbol e_0=\boldsymbol0$. If
agent $m$ is in state $i$ at time $t$, we represent its state by
$\boldsymbol W_m(t)=\boldsymbol e_i$. Define the reduced simplex
\begin{align}
{\cal D}_{n-1}
&=
\left\{
\boldsymbol w\in\mathbb R_+^{n-1}:
\sum_{i=1}^{n-1}w_i\le1
\right\}.
\label{eq:reduced-simplex}
\end{align}
Thus,
$\boldsymbol W_m(t)\in
\{\boldsymbol e_0,\ldots,\boldsymbol e_{n-1}\}
\subset{\cal D}_{n-1}$.
Let
\begin{align}
{\cal W}
&\subseteq
\{\boldsymbol e_0,\ldots,\boldsymbol e_{n-1}\}^M
\label{eq:finite-state-space}
\end{align}
denote the finite {\em admissible} state space of the CTMC, which  allows both unconstrained product state spaces and constrained state spaces. For
$\boldsymbol W\in{\cal W}$, we stack it to a $(n-1)M$ vector such that 
\begin{align}
\boldsymbol W
&=
\left(
\boldsymbol W_1^\top,\ldots,\boldsymbol W_M^\top
\right)^\top
\in\mathbb R^{(n-1)M}.
\label{eq:full-state}
\end{align}
For a bounded weight function $\varphi:\Theta\to\mathbb R$, define the
weighted empirical occupancy
\begin{align}
\boldsymbol X_\varphi(\boldsymbol W)
&=
\frac1M\sum_{m=1}^M
\varphi(\theta_m)\boldsymbol W_m.
\label{eq:weighted-occupancy}
\end{align}
The interaction field is
\begin{align}
\boldsymbol Z(\boldsymbol W)
&=
\Psi\!\left(
\boldsymbol X_\chi(\boldsymbol W)
\right),
\label{eq:interaction-field}
\end{align}
where $\chi:\Theta\to\mathbb R$ is a fixed bounded interaction weight satisfying $0<|\chi|_\infty<\infty,$ and
$\Psi:\mathbb R^{n-1}\to\mathbb R^d$.

We assume that $\boldsymbol W(t)$ is an irreducible CTMC on ${\cal W}$ with
generator ${\cal G}$. Its dynamics are described through a collection
${\cal R}$ of transition types, which may depend on $M$ but not on the
current system state. For each $r\in{\cal R}$, let
$q_r^{\rm CTMC}(\boldsymbol W)$ denote the transition rate, with
$q_r^{\rm CTMC}(\boldsymbol W)=0$ whenever transition $r$ is infeasible at
$\boldsymbol W$, and let
$\boldsymbol\Delta_r\in\mathbb R^{(n-1)M}$ denote the corresponding state
change, which is independent of $\boldsymbol W$. For a test function $g$,
\begin{align}
({\cal G}g)(\boldsymbol W)
&=
\sum_{r\in{\cal R}}
q_r^{\rm CTMC}(\boldsymbol W)
\left(
g(\boldsymbol W+\boldsymbol\Delta_r)-g(\boldsymbol W)
\right).
\label{eq:generator}
\end{align}
The infinitesimal drift and instantaneous covariance-rate matrix of the CTMC
are
\begin{align}
\boldsymbol F^{\rm CTMC}(\boldsymbol W)
&=
\sum_{r\in{\cal R}}
q_r^{\rm CTMC}(\boldsymbol W)\boldsymbol\Delta_r,
\qquad
\boldsymbol Q^{\rm CTMC}(\boldsymbol W)
=
\sum_{r\in{\cal R}}
q_r^{\rm CTMC}(\boldsymbol W)
\boldsymbol\Delta_r\boldsymbol\Delta_r^\top.
\label{eq:ctmc-drift-noise}
\end{align}
Thus, conditional on $\boldsymbol W(t)=\boldsymbol W$,
\begin{align}
\mathbb E[d\boldsymbol W(t)\mid\boldsymbol W(t)=\boldsymbol W]
&=
\boldsymbol F^{\rm CTMC}(\boldsymbol W)\,dt+o(dt),
\qquad
\mathbb E[
d\boldsymbol W(t)d\boldsymbol W(t)^\top
\mid\boldsymbol W(t)=\boldsymbol W]
=
\boldsymbol Q^{\rm CTMC}(\boldsymbol W)\,dt+o(dt).
\label{eq:infinitesimal-moments}
\end{align}
This paper considers a special class of CTMCs in which each agent interacts with other agents {\em only} through  the interaction field ${\bf Z}(\bf W).$   In particular, for each agent, the dependence of its infinitesimal drift on the other agents is through
the interaction field and for every $\boldsymbol{W}\in {\cal W}$ and $m=1,\cdots, M,$ the $m$th block of the CTMC drift satisfies 
\begin{align}
\left[\boldsymbol F^{\rm CTMC}(\boldsymbol W)\right]_m=({\cal G}\boldsymbol W_m)(\boldsymbol W)=\boldsymbol f\left(\boldsymbol W_m, \boldsymbol Z(\boldsymbol W);\theta_m\right).
\end{align} 

\subsection{Heterogeneous Mean-Field Model}
\label{subsec:heterogeneous-mf}

For each realized parameter $\theta_m$, let
$\boldsymbol w_m(t)\in{\cal D}_{n-1}$ denote the mean-field state of agent
$m$. It approximates the conditional mean of the reduced one-hot state,
$\boldsymbol w_m(t)\approx
\mathbb E[\boldsymbol W_m(t)\mid\boldsymbol\theta]$. Equivalently, for $1\leq i\leq n-1,$
\begin{align}
w_{m,i}(t)
&\approx
\mathbb P\!\left(
\boldsymbol W_m(t)=\boldsymbol e_i
\mid\boldsymbol\theta
\right),
\qquad
w_{m,0}(t)
=
1-\sum_{i=1}^{n-1}w_{m,i}(t).
\label{eq:mean-field-probability}
\end{align}
For a bounded weight $\varphi$, define
\begin{align}
\boldsymbol x_\varphi(\boldsymbol w)
&=
\frac1M\sum_{m=1}^M
\varphi(\theta_m)\boldsymbol w_m,
\qquad
\boldsymbol x(\boldsymbol w)
=
\frac1M\sum_{m=1}^M\boldsymbol w_m,
\qquad
\boldsymbol z(\boldsymbol w)
=
\Psi\!\left(
\boldsymbol x_\chi(\boldsymbol w)
\right).
\label{eq:mean-field-occupancies}
\end{align}

For each $M$ and realized $\boldsymbol\theta$, let ${\cal D}$ denote a compact
convex admissible mean-field state space satisfying
\begin{align}
{\cal W}
&\subseteq
{\cal D}
\subseteq
{\cal D}_{n-1}^{\times M}.
\label{eq:mean-field-domain}
\end{align}
We assume that ${\cal D}$ is forward invariant under the mean-field ODE
defined below. The set ${\cal D}$ may incorporate model-specific conservation
constraints.

Conditional on $\boldsymbol\theta$, the quenched mean-field system is
\begin{align}
\dot{\boldsymbol w}_m(t)
&=
\boldsymbol f\left(
\boldsymbol w_m(t),
\boldsymbol z(\boldsymbol w(t));
\theta_m
\right),
\qquad
m=1,\ldots,M,
\qquad
\dot{\boldsymbol w}(t)
=
\boldsymbol F(\boldsymbol w(t)),
\label{eq:quenched-mf}
\end{align} where the mean-field drift $\boldsymbol F(\boldsymbol w(t))$ is defined by its blocks  $\boldsymbol F_m(\boldsymbol w(t))=\boldsymbol f\left(
\boldsymbol w_m(t),
\boldsymbol z(\boldsymbol w(t));
\theta_m
\right).$ Note that $\boldsymbol F(\boldsymbol W)=\boldsymbol F^{\rm CTMC}(\boldsymbol W)$ for every $\boldsymbol W\in{\cal W}$ by construction. 

For an initial condition $\widehat{\boldsymbol w}\in{\cal D}$, let
$\boldsymbol w(t;\widehat{\boldsymbol w})$ denote the solution of
\eqref{eq:quenched-mf}, and let
\begin{align}
\boldsymbol w^*
&=
\left(
(\boldsymbol w_1^*)^\top,\ldots,(\boldsymbol w_M^*)^\top
\right)^\top
\in{\cal D}
\label{eq:equilibrium-full}
\end{align}
denote its equilibrium. Define
\begin{align}
\boldsymbol x_\varphi^*
&=
\frac1M\sum_{m=1}^M
\varphi(\theta_m)\boldsymbol w_m^*,
\qquad
\boldsymbol z^*
=
\Psi(\boldsymbol x_\chi^*),
\qquad
\boldsymbol f(
\boldsymbol w_m^*,
\boldsymbol z^*;
\theta_m)
=
\boldsymbol0,
\quad
m=1,\ldots,M.
\label{eq:weighted-equilibrium}
\end{align}
For a full-state vector
$\boldsymbol u=
(\boldsymbol u_1^\top,\ldots,\boldsymbol u_M^\top)^\top$, define
\begin{align}
\|\boldsymbol u\|_M^2
&=
\frac1M\sum_{m=1}^M
\|\boldsymbol u_m\|^2,
\qquad
\|\boldsymbol u\|_{\max}
=
\max_{1\le m\le M}\|\boldsymbol u_m\|.
\label{eq:full-state-norms}
\end{align}
All $O(\cdot)$ and $o(\cdot)$ statements below are uniform over realized
parameter vectors and admissible domains satisfying the stated assumptions.

%% file: wfo.tex
\section{First-Order Approximation}
\label{sec:first-order}

The first-order analysis in \cite{Yin_26} is based on bounded rates and jumps, smoothness, stability, and occupancy robustness. In this paper, we further strengthen the robustness condition (A5) to control both the first and second variations of the full mean-field trajectory in a normalized block-$\ell_1$ norm with respect to the initial condition. This strengthened condition will result in an approximation error of $O(M^{-1})$ in Proposition \ref{prop:weighted-first-order} for any bounded weight $\varphi$ with constants that may depend on $\varphi.$

Let
\begin{align}
{\cal K}_{\chi}
&=
\left\{
\boldsymbol x\in\mathbb R^{n-1}:\|\boldsymbol x\|_1\le \|\chi\|_\infty
\right\},
\qquad
\|\boldsymbol u\|_{1,M}
=
\frac1M\sum_{m=1}^M\|\boldsymbol u_m\|,
\label{eq:first-order-Kchi-blocknorm}
\end{align}
where
$\boldsymbol u=(\boldsymbol u_1^\top,\ldots,\boldsymbol u_M^\top)^\top$.
The set ${\cal K}_{\chi}$ is compact, independent of $M$ and
$\boldsymbol\theta$, and contains
$\boldsymbol x_\chi(\boldsymbol w)$ for every $\boldsymbol w\in{\cal D}$.

\begin{assumption}[First Order Assumptions]
\label{ass:first-order}
The following conditions hold.

\begin{itemize}

\item \textbf{A1 -- Bounded transition-rate condition.}
There exists $c>0$, independent of $M$ and the realized parameter vector
$\boldsymbol\theta$, such that
\begin{align}
\sum_{r\in{\cal R}}q_r^{\rm CTMC}(\boldsymbol W)
&\le cM
\label{eq:A1-total-rate}
\end{align}
uniformly over all $\boldsymbol W\in{\cal W}$.

\item \textbf{A2 -- Bounded jump condition.}
There exists $c_e>0$, independent of $M$ and realized parameter vector $\boldsymbol \theta,$ such that each transition changes
the states of at most $c_e$ agents.

\item \textbf{A3 -- Smoothness condition.}
The maps $\Psi$ and $\boldsymbol f(\cdot,\cdot;\theta)$ are twice
continuously differentiable in their state arguments. Their first and
second derivatives are uniformly bounded on ${\cal K}_{\chi}$ and
${\cal D}_{n-1}\times\Psi({\cal K}_{\chi})$, respectively, uniformly over
$\theta\in\Theta$.

\item \textbf{A4 -- Uniform exponential stability.}
There exist $C,\gamma>0$, independent of $M$ and $\boldsymbol\theta$, such
that
\begin{align}
\|\boldsymbol w(t;\widehat{\boldsymbol w})-\boldsymbol w^*\|_M
&\le
Ce^{-\gamma t}
\|\widehat{\boldsymbol w}-\boldsymbol w^*\|_M,
\qquad
\widehat{\boldsymbol w}\in{\cal D},
\quad t\ge0.
\label{eq:A4-exp-stability}
\end{align}

\item \textbf{A5 -- Local perturbation robustness.}
There exist constants $C_1,C_2<\infty$, independent of $M$,
$\boldsymbol\theta$, the system state, the transition type $r$, and
$s\in[0,1]$, such that for every feasible CTMC transition
$\boldsymbol W\to\boldsymbol W+\boldsymbol\Delta_r$ and every
\begin{align}
\widehat{\boldsymbol w}
&=
\boldsymbol W+s\boldsymbol\Delta_r,
\qquad
s\in[0,1],
\label{eq:A5-jump-segment}
\end{align}
the first and second flow variations in direction
$\boldsymbol\Delta_r$ satisfy\footnote{At the endpoints of the segment, the derivatives are understood as the
corresponding one-sided derivatives. Since ${\cal D}$ is convex, the entire
segment in \eqref{eq:A5-jump-segment} lies in ${\cal D}$.}
\begin{align}
\int_0^\infty
\left\|
D_{\widehat{\boldsymbol w}}
\boldsymbol w(t;\widehat{\boldsymbol w})
[\boldsymbol\Delta_r]
\right\|_{1,M}^2dt
&\le
\frac{C_1}{M^2},
\qquad
\int_0^\infty
\left\|
D^2_{\widehat{\boldsymbol w}}
\boldsymbol w(t;\widehat{\boldsymbol w})
[\boldsymbol\Delta_r,\boldsymbol\Delta_r]
\right\|_{1,M}dt
\le
\frac{C_2}{M^2}.
\label{eq:A5-robustness}
\end{align}

\end{itemize}
\end{assumption}

\begin{proposition}[First-order approximation error of a weighted occupancy]
\label{prop:weighted-first-order}
Assume A1--A5 and fix a bounded observation weight
$\varphi:\Theta\to\mathbb R$. Recall $\boldsymbol X_\varphi$ and
$\boldsymbol x_\varphi^*$ from Section~\ref{sec:models}. Let
$h:\mathbb R^{n-1}\to\mathbb R$ be twice continuously differentiable on a
neighborhood of
\begin{align}
{\cal K}_\varphi
&=
\left\{
\boldsymbol y\in\mathbb R^{n-1}:
\|\boldsymbol y\|_1\le\|\varphi\|_\infty
\right\}.
\label{eq:Kphi-first-order}
\end{align}
Then
\begin{align}
\left|
\mathbb E\!\left[
h(\boldsymbol X_\varphi(\infty))
\mid\boldsymbol\theta
\right]
-h(\boldsymbol x_\varphi^*)
\right|
&=
O(M^{-1}),
\label{eq:weighted-first-order-result}
\end{align}
where the constant may depend on $h$ and $\varphi$ but is independent of
$M$ and $\boldsymbol\theta$.
\end{proposition}

\begin{proof}
Condition on $\boldsymbol\theta$ and suppress this conditioning from the
notation. Write
$\boldsymbol x_\varphi(t;\widehat{\boldsymbol w})
=
\boldsymbol x_\varphi(\boldsymbol w(t;\widehat{\boldsymbol w}))$.
For every
$\boldsymbol v=
(\boldsymbol v_1^\top,\ldots,\boldsymbol v_M^\top)^\top
\in{\cal D}$,
\begin{align}
\left\|
\frac1M\sum_{m=1}^M
\varphi(\theta_m)\boldsymbol v_m
\right\|_1
&\le
\frac1M\sum_{m=1}^M
|\varphi(\theta_m)|\|\boldsymbol v_m\|_1
\le
\|\varphi\|_\infty.
\label{eq:weighted-occupancy-compact}
\end{align}
Hence $\boldsymbol X_\varphi(\boldsymbol W)$,
$\boldsymbol x_\varphi(t;\widehat{\boldsymbol w})$, and
$\boldsymbol x_\varphi^*$ belong to ${\cal K}_\varphi$. Define the following two constants based on the operator norm induced by the Euclidean norm:
\begin{align}
B_{1,h}
&=
\sup_{\boldsymbol y\in{\cal K}_\varphi}
\|Dh(\boldsymbol y)\|_{\rm op}=\sup_{\boldsymbol y\in{\cal K}_\varphi}\|\nabla h(\boldsymbol y)\|_2,
\qquad
B_{2,h}
=
\sup_{\boldsymbol y\in{\cal K}_\varphi}
\|D^2h(\boldsymbol y)\|_{\rm op}.
\label{eq:h-derivative-bounds}
\end{align}
These constants are finite.

Define the infinite-horizon Stein function
\begin{align}
g(\widehat{\boldsymbol w})
&=
-\int_0^\infty
\left[
h\!\left(
\boldsymbol x_\varphi(t;\widehat{\boldsymbol w})
\right)
-h(\boldsymbol x_\varphi^*)
\right]dt.
\label{eq:weighted-Stein-function}
\end{align}
By the definition of $\boldsymbol x_\varphi$, Cauchy--Schwarz, and A4,
\begin{align}
\|
\boldsymbol x_\varphi(t;\widehat{\boldsymbol w})
-
\boldsymbol x_\varphi^*
\|
&\le
\|\varphi\|_\infty
\|
\boldsymbol w(t;\widehat{\boldsymbol w})
-
\boldsymbol w^*
\|_M
\le
C\|\varphi\|_\infty e^{-\gamma t}
\|
\widehat{\boldsymbol w}
-
\boldsymbol w^*
\|_M.
\label{eq:weighted-flow-stability}
\end{align}
Thus \eqref{eq:weighted-Stein-function} is finite and $g$ is bounded on
${\cal D}$. The flow property gives
\begin{align}
g\!\left(
\boldsymbol w(s;\widehat{\boldsymbol w})
\right)
&=
-\int_s^\infty
\left[
h\!\left(
\boldsymbol x_\varphi(t;\widehat{\boldsymbol w})
\right)
-h(\boldsymbol x_\varphi^*)
\right]dt,
\label{eq:Stein-flow-identity}
\end{align}
and differentiating at $s=0$ yields the Poisson equation
\begin{align}
Dg(\widehat{\boldsymbol w})
[\boldsymbol F(\widehat{\boldsymbol w})]
&=
h\!\left(
\boldsymbol x_\varphi(0;\widehat{\boldsymbol w})
\right)
-
h(\boldsymbol x_\varphi^*).
\label{eq:weighted-Poisson-equation}
\end{align}
At a discrete CTMC state $\boldsymbol W$,
$\boldsymbol x_\varphi(0;\boldsymbol W)
=\frac1M\sum_{m=1}^M
\varphi(\theta_m)\boldsymbol W_m=
\boldsymbol X_\varphi(\boldsymbol W)$.
Since $\boldsymbol W(\infty)$ is stationary and $g$ is bounded,
\begin{align}
\mathbb E\!\left[
({\cal G}g)(\boldsymbol W(\infty))
\right]
&=
0,
\qquad
\mathbb E\!\left[
h(\boldsymbol X_\varphi(\infty))
-
h(\boldsymbol x_\varphi^*)
\right]
=
\mathbb E\!\left[
Dg(\boldsymbol W(\infty))
[\boldsymbol F(\boldsymbol W(\infty))]
-
({\cal G}g)(\boldsymbol W(\infty))
\right].
\label{eq:weighted-Stein-identity}
\end{align}

Fix a feasible CTMC transition
$\boldsymbol W\to\boldsymbol W+\boldsymbol\Delta_r$ and a point
$\widetilde{\boldsymbol w}
=
\boldsymbol W+s\boldsymbol\Delta_r$ on its segment. Write
$\boldsymbol v(t)
=
D_{\widetilde{\boldsymbol w}}
\boldsymbol w(t;\widetilde{\boldsymbol w})
[\boldsymbol\Delta_r]$
and
$\boldsymbol r(t)
=
D^2_{\widetilde{\boldsymbol w}}
\boldsymbol w(t;\widetilde{\boldsymbol w})
[\boldsymbol\Delta_r,\boldsymbol\Delta_r]$.
Differentiating the finite sum defining $\boldsymbol x_\varphi$ gives
\begin{align}
D_{\widetilde{\boldsymbol w}}
\boldsymbol x_\varphi(t;\widetilde{\boldsymbol w})
[\boldsymbol\Delta_r]
&=
\frac1M
\sum_{m=1}^M
\varphi(\theta_m)\boldsymbol v_m(t),
\qquad
D^2_{\widetilde{\boldsymbol w}}
\boldsymbol x_\varphi(t;\widetilde{\boldsymbol w})
[\boldsymbol\Delta_r,\boldsymbol\Delta_r]
=
\frac1M
\sum_{m=1}^M
\varphi(\theta_m)\boldsymbol r_m(t).
\label{eq:weighted-flow-variations}
\end{align}
Therefore, by the triangle inequality, 
\begin{align}
\left\|
D_{\widetilde{\boldsymbol w}}
\boldsymbol x_\varphi(t;\widetilde{\boldsymbol w})
[\boldsymbol\Delta_r]
\right\|
&\le
\|\varphi\|_\infty
\|\boldsymbol v(t)\|_{1,M},
\qquad
\left\|
D^2_{\widetilde{\boldsymbol w}}
\boldsymbol x_\varphi(t;\widetilde{\boldsymbol w})
[\boldsymbol\Delta_r,\boldsymbol\Delta_r]
\right\|
\le
\|\varphi\|_\infty
\|\boldsymbol r(t)\|_{1,M}.
\label{eq:weighted-flow-variation-bounds}
\end{align}
Squaring the first inequality, integrating, and applying A5 yield
\begin{align}
\int_0^\infty
\left\|
D_{\widetilde{\boldsymbol w}}
\boldsymbol x_\varphi(t;\widetilde{\boldsymbol w})
[\boldsymbol\Delta_r]
\right\|^2dt
&\le
\frac{C_1\|\varphi\|_\infty^2}{M^2},
\qquad
\int_0^\infty
\left\|
D^2_{\widetilde{\boldsymbol w}}
\boldsymbol x_\varphi(t;\widetilde{\boldsymbol w})
[\boldsymbol\Delta_r,\boldsymbol\Delta_r]
\right\|dt
\le
\frac{C_2\|\varphi\|_\infty}{M^2}.
\label{eq:A5-weighted-consequences}
\end{align}

Let $\boldsymbol W$ be a discrete CTMC state and let
$\boldsymbol\Delta_r$ be a feasible transition increment at
$\boldsymbol W$. Since ${\cal D}$ is convex and both
$\boldsymbol W$ and $\boldsymbol W+\boldsymbol\Delta_r$ belong to
${\cal W}\subseteq{\cal D}$, the segment between them remains in
${\cal D}$. Taylor's theorem gives
\begin{align}
g(\boldsymbol W+\boldsymbol\Delta_r)
&=
g(\boldsymbol W)
+
Dg(\boldsymbol W)[\boldsymbol\Delta_r]
+
R(\boldsymbol W,\boldsymbol\Delta_r),
\qquad
R(\boldsymbol W,\boldsymbol\Delta_r)
=
\int_0^1
(1-s)
D^2g(\boldsymbol W+s\boldsymbol\Delta_r)
[\boldsymbol\Delta_r,\boldsymbol\Delta_r]ds.
\label{eq:g-Taylor}
\end{align}
By A3, the first and second flow variations exist, and the integrability
bounds in A5 justify differentiation under the integral in
\eqref{eq:weighted-Stein-function}. Hence
\begin{align}
D^2g(\widetilde{\boldsymbol w})
[\boldsymbol\Delta_r,\boldsymbol\Delta_r]
={}&
-\int_0^\infty
D^2h\!\left(
\boldsymbol x_\varphi(t;\widetilde{\boldsymbol w})
\right)
\left[
D_{\widetilde{\boldsymbol w}}
\boldsymbol x_\varphi(t;\widetilde{\boldsymbol w})
[\boldsymbol\Delta_r],
D_{\widetilde{\boldsymbol w}}
\boldsymbol x_\varphi(t;\widetilde{\boldsymbol w})
[\boldsymbol\Delta_r]
\right]dt
\nonumber\\
&-
\int_0^\infty
Dh\!\left(
\boldsymbol x_\varphi(t;\widetilde{\boldsymbol w})
\right)
\left[
D^2_{\widetilde{\boldsymbol w}}
\boldsymbol x_\varphi(t;\widetilde{\boldsymbol w})
[\boldsymbol\Delta_r,\boldsymbol\Delta_r]
\right]dt.
\label{eq:g-second-derivative}
\end{align}
Using \eqref{eq:A5-weighted-consequences},
\begin{align}
\left|
D^2g(\widetilde{\boldsymbol w})
[\boldsymbol\Delta_r,\boldsymbol\Delta_r]
\right|
&\le
\frac{
B_{2,h}C_1\|\varphi\|_\infty^2
+
B_{1,h}C_2\|\varphi\|_\infty
}{M^2}
=
\frac{C_{h,\varphi}}{M^2}.
\label{eq:g-second-bound}
\end{align}
Substituting the bound above into \eqref{eq:g-Taylor} yields
$|R(\boldsymbol W,\boldsymbol\Delta_r)|
\le C_{h,\varphi}/(2M^2)$ uniformly over feasible transitions. Note that
at every discrete CTMC state $\boldsymbol W$, we have
\begin{align}
({\cal G}g)(\boldsymbol W)
&=
\sum_{r\in{\cal R}}
q_r^{\rm CTMC}(\boldsymbol W)
\left(
g(\boldsymbol W+\boldsymbol\Delta_r)-g(\boldsymbol W)
\right)=
Dg(\boldsymbol W)[\boldsymbol F(\boldsymbol W)]
+
\sum_{r\in{\cal R}}
q_r^{\rm CTMC}(\boldsymbol W)
R(\boldsymbol W,\boldsymbol\Delta_r),
\end{align} which implies that 
\begin{align}
\left|
({\cal G}g)(\boldsymbol W)
-
Dg(\boldsymbol W)[\boldsymbol F(\boldsymbol W)]
\right|
\le
\frac{cC_{h,\varphi}}{2M},
\label{eq:generator-Taylor-bound}
\end{align}
where the inequality follows directly from A1 and \eqref{eq:g-second-bound}. Substituting
\eqref{eq:generator-Taylor-bound} into \eqref{eq:weighted-Stein-identity}
and restoring the conditioning on $\boldsymbol\theta$ proves
\eqref{eq:weighted-first-order-result}.
\end{proof}

\begin{corollary}[Weighted occupancy second-moment and bias bounds]
\label{cor:weighted-first-order-moments}
Under A1--A5, for every fixed bounded observation weight $\varphi$,
\begin{align}
\mathbb E\!\left[
\|
\boldsymbol X_\varphi(\infty)
-
\boldsymbol x_\varphi^*
\|^2
\mid\boldsymbol\theta
\right]
&=
O(M^{-1}),
\qquad
\left\|
\mathbb E\!\left[
\boldsymbol X_\varphi(\infty)
-
\boldsymbol x_\varphi^*
\mid\boldsymbol\theta
\right]
\right\|
=
O(M^{-1}).
\label{eq:weighted-second-bias}
\end{align}
Equivalently, with
$\boldsymbol Y_\varphi = \sqrt M
\left( \boldsymbol X_\varphi(\infty) - \boldsymbol x_\varphi^*\right),$
we have
\begin{align}
\mathbb E\!\left[
\|
\boldsymbol Y_\varphi\|^2
\mid\boldsymbol\theta
\right]
&=
O(1),
\qquad
\left\|
\mathbb E[
\boldsymbol Y_\varphi
\mid\boldsymbol\theta]
\right\|
=
O(M^{-1/2}).
\label{eq:Yphi-bias-first-order}
\end{align}
\end{corollary}

\begin{proof}
For the first statement in \eqref{eq:weighted-second-bias}, apply
Proposition~\ref{prop:weighted-first-order} to
$h_{\boldsymbol\theta}(\boldsymbol y)
=
\|\boldsymbol y-\boldsymbol x_\varphi^*\|^2$.
Its first and second derivatives are uniformly bounded on
${\cal K}_\varphi$, independently of $M$ and $\boldsymbol\theta$.
For the bias statement, apply Proposition~\ref{prop:weighted-first-order}
to the coordinate functions $h_i(\boldsymbol y)=y_i$,
$i=1,\ldots,n-1$. Since $n$ is fixed, the resulting coordinatewise
$O(M^{-1})$ bounds imply the second statement in
\eqref{eq:weighted-second-bias}. Equation
\eqref{eq:Yphi-bias-first-order} follows by the definition of $\boldsymbol Y_\varphi$.
\end{proof}

%% file: so.tex
\section{Second-Order Covariance Approximation and Mean Refinement}
\label{sec:second-order}

For the second-order analysis, let $\boldsymbol Q$ be a symmetric positive
semidefinite matrix-valued function on ${\cal D}$ satisfying
$\boldsymbol Q(\boldsymbol W)=\boldsymbol Q^{\mathrm{CTMC}}(\boldsymbol W)$
for every $\boldsymbol W\in{\cal W}$, where $\boldsymbol Q^{\mathrm{CTMC}}$
is defined in~\eqref{eq:ctmc-drift-noise}. Define
\begin{align}
\boldsymbol J
&=D\boldsymbol F(\boldsymbol w^*),
&
\overline{\boldsymbol Q}
&=\boldsymbol Q(\boldsymbol w^*),
&
\boldsymbol R_F(\boldsymbol w)
&=\boldsymbol F(\boldsymbol w)
-\boldsymbol J(\boldsymbol w-\boldsymbol w^*).
\label{eq:second-order-linearization}
\end{align}
For a deterministic block vector
$\boldsymbol a=(\boldsymbol a_1^\top,\ldots,\boldsymbol a_M^\top)^\top$, define the scaled vector 
\begin{align*}
L_{\boldsymbol a}
=\frac{1}{\sqrt M}
\begin{pmatrix}
\boldsymbol a_1^\top&\cdots&\boldsymbol a_M^\top
\end{pmatrix}.
\end{align*}
For a fixed bounded observation weight $\varphi:\Theta\to\mathbb R$, define
\begin{align}
L_\varphi
&=\frac{1}{\sqrt M}
\begin{pmatrix}
\varphi(\theta_1)\boldsymbol I_{n-1}&\cdots&
\varphi(\theta_M)\boldsymbol I_{n-1}
\end{pmatrix},
&
\boldsymbol Y_\varphi
&=L_\varphi(\boldsymbol W(\infty)-\boldsymbol w^*)
=\sqrt M\bigl(
\boldsymbol X_\varphi(\infty)-\boldsymbol x_\varphi^*
\bigr).
\label{eq:second-order-Yphi}
\end{align}
Whenever $\boldsymbol J$ is Hurwitz, let
$\overline{\boldsymbol\Sigma}$ denote the unique symmetric solution of the Lyapunov equation
\begin{align}
\boldsymbol J\overline{\boldsymbol\Sigma}
+\overline{\boldsymbol\Sigma}\boldsymbol J^\top
+\overline{\boldsymbol Q}=\boldsymbol 0
\end{align} such that 
\begin{align}
\overline{\boldsymbol\Sigma}
&=\int_0^\infty
e^{\boldsymbol Jt}\overline{\boldsymbol Q}
e^{\boldsymbol J^\top t}\,dt.
\label{eq:second-order-lyapunov}
\end{align} We further define
\begin{align}
\boldsymbol\Gamma_\varphi
=L_\varphi\overline{\boldsymbol\Sigma}L_\varphi^\top,\quad 
\boldsymbol U
=\boldsymbol W(\infty)-\boldsymbol w^*,\quad
\boldsymbol\Sigma
=\mathbb E[
\boldsymbol U\boldsymbol U^\top
\mid\boldsymbol\theta].
\label{eq:second-order-Sigma}
\end{align}
For a symmetric matrix $\boldsymbol A\in\mathbb R^{N\times N}$, where $N=(n-1)M,$ we define the {\em Hessian contraction} of the mean-field drift at equilibrium by
\begin{align}
[H(\boldsymbol A)]_i
&=\sum_{a,b=1}^N
\frac{\partial^2F_i}
{\partial w_a\partial w_b}(\boldsymbol w^*)A_{ab},
\qquad i=1,\ldots,N.
\label{eq:second-order-H}
\end{align}
Equivalently, for any $\boldsymbol u,\boldsymbol v\in\mathbb R^N$,
\begin{align}
D^2\boldsymbol F(\boldsymbol w^*)[\boldsymbol u,\boldsymbol v]
&=
H\!\left(
\frac{
\boldsymbol u\boldsymbol v^\top+
\boldsymbol v\boldsymbol u^\top}{2}
\right).
\label{eq:second-order-H-bilinear}
\end{align}

The second-order results below are under additional assumptions A6--A9.
A6 controls the linearized dynamics, while A7 provides the stationary
concentration used in Lemma~1. A8 quantifies the nonlinear-drift and
state-dependent-noise errors needed for the covariance approximation, and A9
quantifies the quadratic drift structure and post-quadratic remainder needed for
the refined mean. Thus Theorem~\ref{thm:covariance-approximation} uses A1--A8, whereas Theorem~\ref{thm:performance-refinement} uses A1--A9.
The theorem statements do not require a particular structural representation
of the drift or noise. In Section~\ref{subsec:structural-modeling}, we introduce S1--S2 as convenient
structural conditions that, together with A1--A7, imply A8--A9 through
Lemma~\ref{lem:structural-A8-A9}.

\begin{assumption}[Second-order conditions]
\label{ass:second-order-conditions}
The following conditions hold.

\emph{Stability and concentration.}
\begin{itemize}
\item \textbf{A6 -- Uniform exponential stability of the linearized system.}
There exist constants $C_s,\gamma_s>0$, independent of $M$ and $\boldsymbol\theta$, such that
\begin{align}
\|e^{\boldsymbol Jt}\boldsymbol u\|_{1,M}
&\le
C_se^{-\gamma_st}\|\boldsymbol u\|_{1,M},
\qquad
\boldsymbol u\in\mathbb R^N,
\quad t\ge0.
\label{eq:A6-linear-stability}
\end{align}
In particular, $\boldsymbol J$ is Hurwitz.

\item \textbf{A7 -- Uniform local coupling contraction.}
Recall that ${\cal W}$ is the finite state space of the CTMC. Call two distinct states $\boldsymbol w,\widetilde{\boldsymbol w}\in{\cal W}$ \emph{neighboring}, written $\boldsymbol w\sim\widetilde{\boldsymbol w}$, if one can be obtained from the other by one feasible CTMC transition, in either direction. There is a metric $d$ on ${\cal W}$ and constants $D_e,\kappa_c>0$, independent of $M$ and $\boldsymbol\theta$, such that
\begin{align}
1
&\le
d(\boldsymbol w,\widetilde{\boldsymbol w})
\le D_e\hbox{ if }
\boldsymbol w\sim\widetilde{\boldsymbol w}
\label{eq:A7-edge-diameter}
\end{align}
and every pair $\boldsymbol w,\widetilde{\boldsymbol w}\in{\cal W}$ can be joined by a neighboring path $\boldsymbol w=\boldsymbol w^{(0)},\ldots,\boldsymbol w^{(L)}=\widetilde{\boldsymbol w}$ satisfying
\begin{align}
d(\boldsymbol w,\widetilde{\boldsymbol w})
&=
\sum_{\ell=1}^L
d\!\left(\boldsymbol w^{(\ell-1)},\boldsymbol w^{(\ell)}\right).
\label{eq:A7-shortest-path}
\end{align}
For every neighboring pair there exists a coupling of two copies of the finite CTMC, started from $\boldsymbol w\sim \widetilde{\boldsymbol w}$, whose coupled generator ${\cal G}^c$ satisfies
\begin{align}
({\cal G}^cd)(\boldsymbol w,\widetilde{\boldsymbol w})
&\le
-\kappa_cd(\boldsymbol w,\widetilde{\boldsymbol w}).
\label{eq:A7-local-contraction}
\end{align}

\item \textbf{A8 -- Covariance error control.}
Let
\begin{align}
\boldsymbol E_R
&=\mathbb E\!\left[
\boldsymbol R_F(\boldsymbol W(\infty))\boldsymbol U^\top
+\boldsymbol U\boldsymbol R_F(\boldsymbol W(\infty))^\top
\mid\boldsymbol\theta\right],
&
\boldsymbol E_Q
&=\mathbb E\!\left[
\boldsymbol Q(\boldsymbol W(\infty))-\overline{\boldsymbol Q}
\mid\boldsymbol\theta\right].
\label{eq:A8-errors}
\end{align}
Then
\begin{align}
\sup_{\|\boldsymbol b\|_{\max}\leq1}
\frac1M|\boldsymbol b^\top\boldsymbol E_R\boldsymbol b|
&=O(M^{-1/2}),
&
\sup_{\|\boldsymbol b\|_{\max}\leq1}
\frac1M|\boldsymbol b^\top\boldsymbol E_Q\boldsymbol b|
&=O(M^{-1/2}).
\label{eq:A8-control}
\end{align}

\item \textbf{A9 -- Mean-refinement error control.}
There exist constants $K_H,C_H<\infty$, independent of $M$ and
$\boldsymbol\theta$, such that for every $\|\boldsymbol b\|_{\max}\leq1$
and every symmetric matrix $\boldsymbol A$, there are block vectors
$\boldsymbol r_{k,\boldsymbol b},\boldsymbol s_{k,\boldsymbol b}$ satisfying
\begin{align}
\max_{1\leq k\leq K_H}
\left\{
\|\boldsymbol r_{k,\boldsymbol b}\|_{\max},
\|\boldsymbol s_{k,\boldsymbol b}\|_{\max}
\right\}
&\leq C_H,
&
|\boldsymbol b^\top H(\boldsymbol A)|
&\leq
C_H\sum_{k=1}^{K_H}
\left|
L_{\boldsymbol r_{k,\boldsymbol b}}
\boldsymbol A
L_{\boldsymbol s_{k,\boldsymbol b}}^\top
\right|.
\label{eq:A9-H-control}
\end{align}
Moreover, with
\begin{align}
\boldsymbol T_F(\boldsymbol w)
&=\boldsymbol R_F(\boldsymbol w)
-\frac12H\!\left(
(\boldsymbol w-\boldsymbol w^*)
(\boldsymbol w-\boldsymbol w^*)^\top
\right),
\label{eq:A9-TF}
\end{align}
we have
\begin{align}
\sup_{\|\boldsymbol b\|_{\max}\leq1}
\left|
\boldsymbol b^\top
\mathbb E[\boldsymbol T_F(\boldsymbol W(\infty))\mid\boldsymbol\theta]
\right|
&=o(1).
\label{eq:A9-TF-control}
\end{align}
For the sharper mean-refinement rate, assume the left-hand side of
\eqref{eq:A9-TF-control} is $O(M^{-1/2})$.

\end{itemize}

\end{assumption}

A6--A7 are analytical conditions controlling the linearized dynamics and
stationary concentration, while A8--A9 are the second-order error conditions
used directly in Theorems~\ref{thm:covariance-approximation} and \ref{thm:performance-refinement}. Section~\ref{subsec:structural-modeling} gives structural sufficient
conditions S1--S2 for verifying A8--A9. The local coupling condition A7 is
only one way to obtain the concentration bounds below; any alternative
argument yielding the same uniform moment bounds can be substituted for A7.

Recall that $N=M(n-1)$ denotes the dimension of the full-state vector. For a deterministic block vector
$\boldsymbol a=(\boldsymbol a_1^\top,\ldots,\boldsymbol a_M^\top)^\top$,
define
\begin{align}
x_{\boldsymbol a}(\boldsymbol w)
=\frac{1}{M}\boldsymbol a^\top\boldsymbol w, &\qquad x^*_{\boldsymbol a}=x_{\boldsymbol a}(\boldsymbol w^*),&\qquad \bar{x}_{\boldsymbol a}(\boldsymbol w)=x_{\boldsymbol a}(\boldsymbol w)-x^*_{\boldsymbol a}.
\label{eq:second-order-projections}
\end{align}
The coefficient vector $\boldsymbol a$ may depend on the realized
$\boldsymbol\theta$, but is deterministic after conditioning on
$\boldsymbol\theta$.

\begin{lemma}[Uniform block moment bounds]
\label{lem:uniform-block-moments} 
Under A1--A7, there exists $C<\infty$, independent of $M$, $\boldsymbol\theta$, and $\boldsymbol a$, such that
\begin{align}
\sup_{\|\boldsymbol a\|_{\max}\le1}
\left|\mathbb E[\bar{x}_a(\boldsymbol W(\infty))\mid\boldsymbol\theta]\right|
\le
\frac{C}{M},
&\qquad \sup_{\|\boldsymbol a\|_{\max}\le1}
\mathbb E[|\bar{x}_a(\boldsymbol W(\infty))|^2\mid\boldsymbol\theta]
\le
\frac{C}{M},
&\sup_{\|\boldsymbol a\|_{\max}\le1}
\mathbb E[|\bar{x}_a(\boldsymbol W(\infty))|^4\mid\boldsymbol\theta]
\le
\frac{C}{M^2}.
\label{eq:block-moment-bounds}
\end{align}
\end{lemma}

\begin{proof}
Condition on $\boldsymbol\theta$ and suppress this conditioning and consider $\|\boldsymbol a\|_{\max}\le1$.

We first prove the bias bound using A1--A5. The proof is similar to the proof of Proposition \ref{prop:weighted-first-order} but needs to establish a bound uniform over all block vectors a satisfying
the stated block-max norm restriction. 
Define
\begin{align}
x_{\boldsymbol a}(t;\widehat{\boldsymbol w})
&=
\frac1M\sum_{m=1}^M
\boldsymbol a_m^\top\boldsymbol w_m(t;\widehat{\boldsymbol w}),
\qquad
g_{\boldsymbol a}(\widehat{\boldsymbol w})
=
-\int_0^\infty
\left(x_{\boldsymbol a}(t;\widehat{\boldsymbol w})-x_{\boldsymbol a}^*\right)dt.
\label{eq:block-Stein-function-sec4}
\end{align}
By A4 the integral is finite. Arguing as in the proof of
Proposition~\ref{prop:weighted-first-order}, A3 and A5 justify
differentiation under the integral along every feasible jump segment.
Since $x_{\boldsymbol a}$ is linear, for every feasible CTMC transition
$\boldsymbol W\to\boldsymbol W+\boldsymbol\Delta_r$ and every
$\widehat{\boldsymbol w}
=
\boldsymbol W+s\boldsymbol\Delta_r$, $s\in[0,1]$,
\begin{align}
\left|
D^2g_{\boldsymbol a}(\widehat{\boldsymbol w})
[\boldsymbol\Delta_r,\boldsymbol\Delta_r]
\right|
&\le
\int_0^\infty
\left\|
D^2_{\widehat{\boldsymbol w}}
\boldsymbol w(t;\widehat{\boldsymbol w})
[\boldsymbol\Delta_r,\boldsymbol\Delta_r]
\right\|_{1,M}dt
\le
\frac{C_2}{M^2},
\label{eq:block-Stein-second-bound-sec4}
\end{align} where the last inequality holds due to A5. 
The Poisson equation gives
\begin{align}
Dg_{\boldsymbol a}(\widehat{\boldsymbol w})
[\boldsymbol F(\widehat{\boldsymbol w})]
=
x_{\boldsymbol a}(\widehat{\boldsymbol w})-x_a^*=\bar{x}_{\boldsymbol a}(\widehat{\boldsymbol w}).
\label{eq:ga-poisson}
\end{align}
For a feasible transition
$\boldsymbol W\to\boldsymbol W+\boldsymbol\Delta_r$,
Taylor's theorem gives
\begin{align}
g_{\boldsymbol a}(\boldsymbol W+\boldsymbol\Delta_r)
-g_{\boldsymbol a}(\boldsymbol W)
&=
Dg_{\boldsymbol a}(\boldsymbol W)[\boldsymbol\Delta_r]
+R_r(\boldsymbol W),
\label{eq:ga-taylor}
\\
R_r(\boldsymbol W)
&=
\int_0^1(1-s)
D^2g_{\boldsymbol a}(\boldsymbol W+s\boldsymbol\Delta_r)
[\boldsymbol\Delta_r,\boldsymbol\Delta_r]\,ds.
\label{eq:ga-remainder}
\end{align}
By \eqref{eq:ga-remainder} and \eqref{eq:block-Stein-second-bound-sec4},
\begin{align}
|R_r(\boldsymbol W)|
\leq
\frac{C_2}{2M^2}.
\label{eq:ga-remainder-bound}
\end{align}
Hence, at every discrete CTMC state $\boldsymbol W$,
\begin{align}
({\cal G}g_{\boldsymbol a})(\boldsymbol W)
&=
\sum_{r\in{\cal R}}
q_r^{\rm CTMC}(\boldsymbol W)
\left(
Dg_{\boldsymbol a}(\boldsymbol W)[\boldsymbol\Delta_r]
+R_r(\boldsymbol W)
\right)
\nonumber\\
&=
Dg_{\boldsymbol a}(\boldsymbol W)[\boldsymbol F(\boldsymbol W)]
+
\sum_{r\in{\cal R}}
q_r^{\rm CTMC}(\boldsymbol W)R_r(\boldsymbol W),
\label{eq:ga-generator}
\end{align}
where we used
$\boldsymbol F(\boldsymbol W)
=
\sum_{r\in{\cal R}}
q_r^{\rm CTMC}(\boldsymbol W)\boldsymbol\Delta_r$
on the discrete state space. Since $\boldsymbol W(\infty)$ is stationary,
\begin{align}
\mathbb E[
({\cal G}g_{\boldsymbol a})(\boldsymbol W(\infty))
]
=0.
\label{eq:ga-stationarity}
\end{align}
Combining \eqref{eq:ga-poisson}, \eqref{eq:ga-generator}, and
\eqref{eq:ga-stationarity} yields
\begin{align}
\mathbb E[\bar{x}_{\boldsymbol a}(\boldsymbol W(\infty))]
=
-\mathbb E\!\left[
\sum_{r\in{\cal R}}
q_r^{\rm CTMC}(\boldsymbol W(\infty))
R_r(\boldsymbol W(\infty))
\right].
\label{eq:ga-bias-identity}
\end{align}
Therefore, by \eqref{eq:ga-remainder-bound} and A1,
\begin{align}
|\mathbb E[\bar{x}_{\boldsymbol a}(\boldsymbol W(\infty))]|
&\leq
\frac{C_2}{2M^2}
\mathbb E\!\left[
\sum_{r\in{\cal R}}
q_r^{\rm CTMC}(\boldsymbol W(\infty))
\right]
\leq
\frac{cC_2}{2M}.
\label{eq:ga-bias-bound}
\end{align}

We next use A7 to bound the second and fourth moments. 
Adding and subtracting
$\mathbb E[x_{\boldsymbol a}(\boldsymbol W(\infty))]$ therefore gives the
exact decomposition
\begin{align}
\bar{x}_{\boldsymbol a}(\boldsymbol W(\infty))
&=
\underbrace{
x_{\boldsymbol a}(\boldsymbol W(\infty))
-\mathbb E[x_{\boldsymbol a}(\boldsymbol W(\infty))]
}_{\text{centered fluctuation}}
+
\underbrace{
\mathbb E[x_{\boldsymbol a}(\boldsymbol W(\infty))]
-x_{\boldsymbol a}^*
}_{\mathbb E[\bar{\boldsymbol x}_{\boldsymbol a}(\boldsymbol W(\infty))]}.
\label{eq:fa-bias-fluctuation-decomposition}
\end{align}
Hence, for $p\in\{2,4\}$, the inequality
$|x+y|^p\leq2^{p-1}(|x|^p+|y|^p)$ yields
\begin{align}
\mathbb E|\bar{x}_{\boldsymbol a}(\boldsymbol W(\infty))|^p
&\leq
2^{p-1}\mathbb E\!\left[
\left|
x_{\boldsymbol a}(\boldsymbol W(\infty))
-\mathbb E[x_{\boldsymbol a}(\boldsymbol W(\infty))]
\right|^p
\right]
+2^{p-1}|\mathbb E[\bar{x}_{\boldsymbol a}(\boldsymbol W(\infty))]|^p.
\label{eq:fa-bias-fluctuation-moment-decomposition}
\end{align}
By \eqref{eq:ga-bias-bound}, the second term in
\eqref{eq:fa-bias-fluctuation-moment-decomposition} is $O(M^{-2})$ for
$p=2$ and $O(M^{-4})$ for $p=4$, both smaller than the corresponding
bounds required in the lemma. Therefore it remains only to prove
\begin{align}
\mathbb E\!\left[
\left|
x_{\boldsymbol a}(\boldsymbol W(\infty))
-\mathbb E[x_{\boldsymbol a}(\boldsymbol W(\infty))]
\right|^2
\right]
&\leq
\frac{C}{M},
&
\mathbb E\!\left[
\left|
x_{\boldsymbol a}(\boldsymbol W(\infty))
-\mathbb E[x_{\boldsymbol a}(\boldsymbol W(\infty))]
\right|^4
\right]
&\leq
\frac{C}{M^2}.
\label{eq:fa-centered-goal}
\end{align}

Let $P_t=e^{t{\cal G}}$ denote the time-$t$ transition matrix of the CTMC.
For any function $f:{\cal W}\to\mathbb R$,
\begin{align}
(P_tf)(\boldsymbol w)
&=
\mathbb E\!\left[
f(\boldsymbol W(t))
\mid
\boldsymbol W(0)=\boldsymbol w
\right],
\qquad
P_{s+t}=P_sP_t.
\label{eq:transition-semigroup}
\end{align}

Fix $T>0$ and start the chain in stationarity. Then
$\boldsymbol W(T)$ has the same distribution as $\boldsymbol W(\infty)$,
and
\begin{align}
x_{\boldsymbol a}(\boldsymbol W(T))
-\mathbb E[x_{\boldsymbol a}(\boldsymbol W(\infty))]
&=
\underbrace{
x_{\boldsymbol a}(\boldsymbol W(T))
-P_Tx_{\boldsymbol a}(\boldsymbol W(0))
}_{\text{prediction error}}+
\underbrace{
P_Tx_{\boldsymbol a}(\boldsymbol W(0))
-\mathbb E[x_{\boldsymbol a}(\boldsymbol W(\infty))]
}_{\text{mixing error}}.
\label{eq:fa-stationary-decomposition}
\end{align}
Hence, for $p\in\{2,4\}$ and any $T>0,$
\begin{align}
&\mathbb E\!\left[
\left|
x_{\boldsymbol a}(\boldsymbol W(\infty))
-\mathbb E[x_{\boldsymbol a}(\boldsymbol W(\infty))]
\right|^p
\right]\nonumber\\
\leq&
2^{p-1}\mathbb E\!\left[
\left|
x_{\boldsymbol a}(\boldsymbol W(T))
-P_Tx_{\boldsymbol a}(\boldsymbol W(0))
\right|^p
\right]
+
2^{p-1}\sup_{\boldsymbol w\in{\cal W}}
\left|
P_Tx_{\boldsymbol a}(\boldsymbol w)
-\mathbb E[x_{\boldsymbol a}(\boldsymbol W(\infty))]
\right|^p.
\label{eq:fa-centered-prelimit}
\end{align}
We control the two terms on the right-hand side of
\eqref{eq:fa-centered-prelimit} separately. We first use A7 to show that
the mixing error vanishes exponentially in $T$. We then bound the second
and fourth moments of the prediction error uniformly in $T$ and the
initial state.

We begin with the mixing error. For $t\geq0$, define
\begin{align}
K_{\boldsymbol a}(t)
&=
\sup_{\boldsymbol w\sim\widetilde{\boldsymbol w}}
\frac{
|P_tx_{\boldsymbol a}(\boldsymbol w)
-P_tx_{\boldsymbol a}(\widetilde{\boldsymbol w})|
}{
d(\boldsymbol w,\widetilde{\boldsymbol w})
}.
\label{eq:fa-Lt}
\end{align}
We first extend this neighboring-pair bound to arbitrary pairs. Fix
$\boldsymbol w,\widetilde{\boldsymbol w}\in{\cal W}$. By A7, there exists
a neighboring path
$\boldsymbol w=\boldsymbol w^{(0)},\ldots,
\boldsymbol w^{(L)}=\widetilde{\boldsymbol w}$ such that
\begin{align}
d(\boldsymbol w,\widetilde{\boldsymbol w})
&=
\sum_{\ell=1}^L
d(\boldsymbol w^{(\ell-1)},\boldsymbol w^{(\ell)}).
\label{eq:fa-shortest-path}
\end{align}
For each edge on this path, the definition of $K_{\boldsymbol a}(t)$ gives
\begin{align}
|P_tx_{\boldsymbol a}(\boldsymbol w^{(\ell-1)})
-P_tx_{\boldsymbol a}(\boldsymbol w^{(\ell)})|
&\leq
K_{\boldsymbol a}(t)
d(\boldsymbol w^{(\ell-1)},\boldsymbol w^{(\ell)}),
\qquad
\ell=1,\ldots,L.
\label{eq:fa-edge-bound}
\end{align}
Therefore, by telescoping, the triangle inequality, and
\eqref{eq:fa-shortest-path},
\begin{align}
|P_tx_{\boldsymbol a}(\boldsymbol w)
-P_tx_{\boldsymbol a}(\widetilde{\boldsymbol w})|
&\leq
K_{\boldsymbol a}(t)d(\boldsymbol w,\widetilde{\boldsymbol w}).
\label{eq:fa-path-bound}
\end{align}

Fix a neighboring pair $\boldsymbol w\sim\widetilde{\boldsymbol w}.$ Let $\boldsymbol W^{\boldsymbol w}$ and $\boldsymbol W^{\widetilde{\boldsymbol w}}$ denote the two copies of the CTMC under the coupling in A7, started from $\boldsymbol w$ and $\widetilde{\boldsymbol w},$ respectively. 
The expansion of the coupled-generator, denoted by ${\cal G}^c,$ gives
\begin{align}
\mathbb E\!\left[
d\!\left(
\boldsymbol W^{\boldsymbol w}(h),
\boldsymbol W^{\widetilde{\boldsymbol w}}(h)
\right)
\right]
&=
d(\boldsymbol w,\widetilde{\boldsymbol w})
+h({\cal G}^cd)(\boldsymbol w,\widetilde{\boldsymbol w})+o(h)
\nonumber\\
&\leq
(1-\kappa_ch)d(\boldsymbol w,\widetilde{\boldsymbol w})+o(h).
\label{eq:fa-coupled-generator-expansion}
\end{align}
For each fixed $M$, there are finitely many neighboring pairs, so the
$o(h)$ term can be chosen uniformly over them. Since
$d(\boldsymbol w,\widetilde{\boldsymbol w})\geq1$ for neighboring pairs,
there exists $\varepsilon(h)\to0$ as $h\downarrow0$ such that
\begin{align}
\mathbb E\!\left[
d\!\left(
\boldsymbol W^{\boldsymbol w}(h),
\boldsymbol W^{\widetilde{\boldsymbol w}}(h)
\right)
\right]
&\leq
\bigl(1-\kappa_ch+h\varepsilon(h)\bigr)
d(\boldsymbol w,\widetilde{\boldsymbol w}).
\label{eq:fa-distance-contraction-small-time}
\end{align}

By the Markov property, $P_{t+h}=P_hP_t$. Since the coupling in A7 has
CTMCs started from $\boldsymbol w$ and $\widetilde{\boldsymbol w}$ as its
two marginals,
\begin{align}
P_{t+h}x_{\boldsymbol a}(\boldsymbol w)
&=
\mathbb E\!\left[
P_tx_{\boldsymbol a}(\boldsymbol W^{\boldsymbol w}(h))
\right],
&
P_{t+h}x_{\boldsymbol a}(\widetilde{\boldsymbol w})
&=
\mathbb E\!\left[
P_tx_{\boldsymbol a}(\boldsymbol W^{\widetilde{\boldsymbol w}}(h))
\right].
\label{eq:fa-markov-identity}
\end{align}
Hence, by \eqref{eq:fa-path-bound} and
\eqref{eq:fa-distance-contraction-small-time}, the following inequality holds for any $\boldsymbol{w}\sim \tilde{\boldsymbol{w}}$
\begin{align}
|P_{t+h}x_{\boldsymbol a}(\boldsymbol w)
-P_{t+h}x_{\boldsymbol a}(\widetilde{\boldsymbol w})|
&\leq
\bigl(1-\kappa_ch+h\varepsilon(h)\bigr)
K_{\boldsymbol a}(t)d(\boldsymbol w,\widetilde{\boldsymbol w}),
\label{eq:fa-L-recursion-step}
\end{align}
and therefore
\begin{align}
K_{\boldsymbol a}(t+h)
&\leq
\bigl(1-\kappa_ch+h\varepsilon(h)\bigr)
K_{\boldsymbol a}(t).
\label{eq:fa-L-recursion}
\end{align}

The explicit linear form of $x_{\boldsymbol a}$ and A2 give the initial
bound. For a neighboring pair $\boldsymbol w\sim\widetilde{\boldsymbol w}$,
define
\begin{align}
{\cal I}(\boldsymbol w,\widetilde{\boldsymbol w})
&=
\{m:\boldsymbol w_m\neq\widetilde{\boldsymbol w}_m\}.
\label{eq:fa-changed-blocks}
\end{align}
Since $P_0x_{\boldsymbol a}=x_{\boldsymbol a}$, the exact initial value is
\begin{align}
K_{\boldsymbol a}(0)
&=
\frac{1}{M}
\sup_{\boldsymbol w\sim\widetilde{\boldsymbol w}}
\frac{
\left|
\sum_{m\in{\cal I}(\boldsymbol w,\widetilde{\boldsymbol w})}
\boldsymbol a_m^\top
(\boldsymbol w_m-\widetilde{\boldsymbol w}_m)
\right|
}{
d(\boldsymbol w,\widetilde{\boldsymbol w})
}.
\label{eq:fa-L0-exact}
\end{align}
By A2,
$|{\cal I}(\boldsymbol w,\widetilde{\boldsymbol w})|\leq c_e$.
Since each block is a reduced one-hot state, any two block states satisfy
$\|\boldsymbol w_m-\widetilde{\boldsymbol w}_m\|\leq\sqrt{2}$.
Therefore, for $\|\boldsymbol a\|_{\max}\leq1$,
\begin{align}
\left|
\sum_{m\in{\cal I}(\boldsymbol w,\widetilde{\boldsymbol w})}
\boldsymbol a_m^\top
(\boldsymbol w_m-\widetilde{\boldsymbol w}_m)
\right|
&\leq
\sum_{m\in{\cal I}(\boldsymbol w,\widetilde{\boldsymbol w})}
\|\boldsymbol a_m\|\,
\|\boldsymbol w_m-\widetilde{\boldsymbol w}_m\|
\leq
\sqrt{2}\,c_e.
\label{eq:fa-neighbor-initial-bound}
\end{align}
Since A7 gives $d(\boldsymbol w,\widetilde{\boldsymbol w})\geq1$ for
neighboring states, \eqref{eq:fa-L0-exact} yields
\begin{align}
K_{\boldsymbol a}(0)
&\leq
\frac{\sqrt{2}\,c_e}{M}.
\label{eq:fa-L0}
\end{align}

Fix $t>0$ and take $h=t/k$. For all sufficiently large $k$, iterating
\eqref{eq:fa-L-recursion} $k$ times gives
\begin{align}
K_{\boldsymbol a}(t)
&\leq
\left(
1-\frac{\kappa_ct}{k}
+\frac{t}{k}\varepsilon(t/k)
\right)^k
K_{\boldsymbol a}(0).
\label{eq:fa-L-iteration}
\end{align}
Since $\varepsilon(t/k)\to0$, the right-hand side of
\eqref{eq:fa-L-iteration} converges to
$e^{-\kappa_ct}K_{\boldsymbol a}(0)$. Since the left-hand side does not
depend on $k$, taking $k\to\infty$ yields
\begin{align}
K_{\boldsymbol a}(t)
&\leq
e^{-\kappa_ct}K_{\boldsymbol a}(0).
\label{eq:fa-L-limit}
\end{align}
Combining \eqref{eq:fa-path-bound}, \eqref{eq:fa-L0}, and
\eqref{eq:fa-L-limit},
\begin{align}
|P_tx_{\boldsymbol a}(\boldsymbol w)
-P_tx_{\boldsymbol a}(\widetilde{\boldsymbol w})|
&\leq
\frac{\sqrt{2}\,c_e}{M}
e^{-\kappa_ct}
d(\boldsymbol w,\widetilde{\boldsymbol w}),
\qquad
t\geq0.
\label{eq:fa-semigroup-contraction}
\end{align}

The contraction in \eqref{eq:fa-semigroup-contraction} now controls the
mixing error in \eqref{eq:fa-centered-prelimit}. By stationarity,
\begin{align}
\mathbb E[P_tx_{\boldsymbol a}(\boldsymbol W(\infty))]
&=
\mathbb E[x_{\boldsymbol a}(\boldsymbol W(\infty))].
\label{eq:fa-stationary-semigroup-mean}
\end{align}
Therefore, by \eqref{eq:fa-semigroup-contraction} and A7,
\begin{align}
\left|
P_tx_{\boldsymbol a}(\boldsymbol w)
-\mathbb E[x_{\boldsymbol a}(\boldsymbol W(\infty))]
\right|
&=
\left|
P_tx_{\boldsymbol a}(\boldsymbol w)
-\mathbb E[P_tx_{\boldsymbol a}(\boldsymbol W(\infty))]
\right|
\nonumber\\
&\leq
\mathbb E\!\left[
\left|
P_tx_{\boldsymbol a}(\boldsymbol w)
-P_tx_{\boldsymbol a}(\boldsymbol W(\infty))
\right|
\right]
\nonumber\\
&\leq
\frac{\sqrt{2}\,c_e}{M}e^{-\kappa_ct}
\mathbb E[
d(\boldsymbol w,\boldsymbol W(\infty))
],
\end{align}
which implies that 
\begin{align}
\sup_{\boldsymbol{w}\in{\cal W}}\left|
P_Tx_{\boldsymbol a}(\boldsymbol w)
-\mathbb E[x_{\boldsymbol a}(\boldsymbol W(\infty))]
\right|\leq
\frac{\sqrt{2}c_e}{M}\hbox{diam}({\cal W})e^{-\kappa_cT}.
\label{eq:fa-stationary-mixing}
\end{align}

By \eqref{eq:fa-stationary-mixing}, the second term in
\eqref{eq:fa-centered-prelimit} vanishes exponentially as $T\to\infty$.
Therefore, it remains to prove, uniformly in $T$ and the initial state,
\begin{align}
\mathbb E\!\left[
\left|
x_{\boldsymbol a}(\boldsymbol W(T))
-P_Tx_{\boldsymbol a}(\boldsymbol W(0))
\right|^2
\right]
&\leq
\frac{C}{M},
&
\mathbb E\!\left[
\left|
x_{\boldsymbol a}(\boldsymbol W(T))
-P_Tx_{\boldsymbol a}(\boldsymbol W(0))
\right|^4
\right]
&\leq
\frac{C}{M^2}.
\label{eq:fa-prediction-error-goal}
\end{align}

To establish \eqref{eq:fa-prediction-error-goal}, fix $T>0$ and an
arbitrary initial state, and let ${\cal F}_t$ be the history of the CTMC
up to time $t$. Define the conditional prediction
\begin{align}
{\cal M}_t
&=
\mathbb E\!\left[
x_{\boldsymbol a}(\boldsymbol W(T))
\mid{\cal F}_t
\right]
=
P_{T-t}x_{\boldsymbol a}(\boldsymbol W(t)),
\qquad
0\leq t\leq T.
\label{eq:fa-conditional-prediction}
\end{align}
The second equality follows from the Markov property. Thus
${\cal M}_t$ is the conditional prediction, at time $t$, of the terminal
value $x_{\boldsymbol a}(\boldsymbol W(T))$. For $0\leq s\leq t\leq T$,
the tower property gives
\begin{align}
\mathbb E[{\cal M}_t\mid{\cal F}_s]
&=
\mathbb E\!\left[
\mathbb E\!\left[
x_{\boldsymbol a}(\boldsymbol W(T))
\mid{\cal F}_t
\right]
\Bigm|{\cal F}_s
\right]
=
\mathbb E\!\left[
x_{\boldsymbol a}(\boldsymbol W(T))
\mid{\cal F}_s
\right]
=
{\cal M}_s.
\label{eq:fa-prediction-martingale}
\end{align}
Hence $({\cal M}_t)_{0\leq t\leq T}$ is a martingale. Define the centered
conditional prediction
\begin{align}
{\cal Z}_t
&=
{\cal M}_t-{\cal M}_0,
\qquad
0\leq t\leq T.
\label{eq:fa-centered-prediction}
\end{align}
Then ${\cal Z}_t$ is a martingale and ${\cal Z}_0=0$.

If a feasible transition of type $r$ occurs from state $\boldsymbol w$ at
time $t$, we further define 
\begin{align}
\delta_r(t,\boldsymbol w)
&=
P_{T-t}x_{\boldsymbol a}(\boldsymbol w+\boldsymbol\Delta_r)
-P_{T-t}x_{\boldsymbol a}(\boldsymbol w).
\label{eq:fa-prediction-jump}
\end{align}
Since the two states are neighboring,
\eqref{eq:fa-semigroup-contraction} and A7 imply
\begin{align}
|\delta_r(t,\boldsymbol w)|
&\leq
\frac{\sqrt{2}\,c_eD_e}{M}e^{-\kappa_c(T-t)}.
\label{eq:fa-prediction-jump-bound}
\end{align}

Since ${\cal Z}_t$ is ${\cal F}_t$-measurable, expanding the square and
conditioning on ${\cal F}_t$ gives
\begin{align}
\mathbb E\!\left[
{\cal Z}_{t+h}^2-{\cal Z}_t^2
\mid{\cal F}_t
\right]
&=
2{\cal Z}_t
\mathbb E\!\left[
{\cal Z}_{t+h}-{\cal Z}_t
\mid{\cal F}_t
\right]
+
\mathbb E\!\left[
({\cal Z}_{t+h}-{\cal Z}_t)^2
\mid{\cal F}_t
\right]
\nonumber\\
&=
\mathbb E\!\left[
({\cal Z}_{t+h}-{\cal Z}_t)^2
\mid{\cal F}_t
\right],
\label{eq:fa-second-increment-expansion}
\end{align}
where the last equality holds because ${\cal Z}_t$ is a martingale. 
For sufficiently small $h>0,$ we have
\begin{align}
\mathbb E\!\left[
{\cal Z}_{t+h}^2-{\cal Z}_t^2
\mid{\cal F}_t
\right]
&=
h\sum_{r\in{\cal R}}
q_r^{\rm CTMC}(\boldsymbol W(t))
\delta_r(t,\boldsymbol W(t))^2
+o(h).
\label{eq:fa-second-increment-rate}
\end{align}
Taking expectations and using the tower property,
\begin{align}
\mathbb E[{\cal Z}_{t+h}^2]-\mathbb E[{\cal Z}_t^2]
&=
h\,\mathbb E\!\left[
\sum_{r\in{\cal R}}
q_r^{\rm CTMC}(\boldsymbol W(t))
\delta_r(t,\boldsymbol W(t))^2
\right]
+o(h)
\label{eq:fa-second-increment-unconditional}
\end{align}
where the remainder remains $o(h)$ after taking expectation because the
remainder in \eqref{eq:fa-second-increment-rate} can be chosen to be uniform over
the finite state space. Dividing by $h$ gives
\begin{align}
\frac{
\mathbb E[{\cal Z}_{t+h}^2]-\mathbb E[{\cal Z}_t^2]
}{h}
&=
\mathbb E\!\left[
\sum_{r\in{\cal R}}
q_r^{\rm CTMC}(\boldsymbol W(t))
\delta_r(t,\boldsymbol W(t))^2
\right]
+o(1).
\label{eq:fa-second-difference-quotient}
\end{align}
Moreover, by A1 and \eqref{eq:fa-prediction-jump-bound}, for every
$\boldsymbol w\in{\cal W}$,
\begin{align}
\sum_{r\in{\cal R}}
q_r^{\rm CTMC}(\boldsymbol w)
\delta_r(t,\boldsymbol w)^2
&\leq
\left(
\sum_{r\in{\cal R}}
q_r^{\rm CTMC}(\boldsymbol w)
\right)
\frac{K^2}{M^2}
e^{-2\kappa_c(T-t)}
\nonumber\\
&\leq
\frac{cK^2}{M}
e^{-2\kappa_c(T-t)}
=:a_2(t),
\label{eq:fa-second-rate-bound}
\end{align} where $K=\sqrt{2}c_e D_e.$
Letting $h\downarrow0$ in
\eqref{eq:fa-second-difference-quotient} and using
\eqref{eq:fa-second-rate-bound},
\begin{align}
\frac{d}{dt}\mathbb E[{\cal Z}_t^2]
&=
\mathbb E\!\left[
\sum_{r\in{\cal R}}
q_r^{\rm CTMC}(\boldsymbol W(t))
\delta_r(t,\boldsymbol W(t))^2
\right]
\leq
a_2(t).
\label{eq:fa-prediction-second-moment}
\end{align}
Since ${\cal Z}_0=0$, 
\begin{align}
\mathbb E[{\cal Z}_t^2]
&\leq
\int_0^t a_2(s)\,ds
\leq
\frac{cK^2}{2\kappa_cM},
\qquad
0\leq t\leq T.
\label{eq:fa-prediction-second-bound}
\end{align}

For the fourth moment, we can similarly obtain 
\begin{align}
\frac{d}{dt}\mathbb E[{\cal Z}_t^4]
&=
\mathbb E\!\left[
\sum_{r\in{\cal R}}
q_r^{\rm CTMC}(\boldsymbol W(t))
\left(
6{\cal Z}_t^2\delta_r^2
+4{\cal Z}_t\delta_r^3
+\delta_r^4
\right)
\right],
\label{eq:fa-prediction-fourth-expansion}
\end{align}
where
$\delta_r=\delta_r(t,\boldsymbol W(t))$.
Since
$4|x||y|^3\leq2x^2y^2+2y^4$,
\begin{align}
\frac{d}{dt}\mathbb E[{\cal Z}_t^4]
&\leq
8a_2(t)\mathbb E[{\cal Z}_t^2]
+3a_4(t),
\qquad
a_4(t)
=
\frac{cK^4}{M^3}
e^{-4\kappa_c(T-t)}.
\label{eq:fa-prediction-fourth-differential}
\end{align}
Integrating \eqref{eq:fa-prediction-fourth-differential} from $0$ to
$t\leq T$ and using \eqref{eq:fa-prediction-second-bound},
\begin{align}
\mathbb E[{\cal Z}_t^4]
&\leq
8\sup_{0\leq s\leq t}
\mathbb E[{\cal Z}_s^2]
\int_0^t a_2(s)\,ds
+
3\int_0^t a_4(s)\,ds
\leq
\frac{C}{M^2},
\qquad
0\leq t\leq T.
\label{eq:fa-prediction-fourth-bound}
\end{align}
Thus
\begin{align}
\mathbb E[{\cal Z}_t^2]
&\leq
\frac{C}{M},
&
\mathbb E[{\cal Z}_t^4]
&\leq
\frac{C}{M^2},
\qquad
0\leq t\leq T.
\label{eq:fa-prediction-moment-summary}
\end{align}

We observe that 
\begin{align}
{\cal Z}_T
&=
x_{\boldsymbol a}(\boldsymbol W(T))
-P_Tx_{\boldsymbol a}(\boldsymbol W(0)).
\label{eq:fa-Z-terminal}
\end{align}
Setting $t=T$ in \eqref{eq:fa-prediction-moment-summary} therefore gives,
uniformly over the initial state,
\begin{align}
\mathbb E[|{\cal Z}_T|^2]
&=
\mathbb E\!\left[
\left|
x_{\boldsymbol a}(\boldsymbol W(T))
-P_Tx_{\boldsymbol a}(\boldsymbol W(0))
\right|^2
\right]
\leq
\frac{C}{M},
&
\mathbb E[|{\cal Z}_T|^4]
&=
\mathbb E\!\left[
\left|
x_{\boldsymbol a}(\boldsymbol W(T))
-P_Tx_{\boldsymbol a}(\boldsymbol W(0))
\right|^4
\right]
\leq
\frac{C}{M^2}.
\label{eq:fa-conditional-moment-bounds}
\end{align}

Returning to the stationary chain, substitute
\eqref{eq:fa-conditional-moment-bounds} into
\eqref{eq:fa-centered-prelimit} and let $T\to\infty$.
The mixing term vanishes by \eqref{eq:fa-stationary-mixing}, so
\begin{align}
\mathbb E\!\left[
\left|
x_{\boldsymbol a}(\boldsymbol W(\infty))
-\mathbb E[x_{\boldsymbol a}(\boldsymbol W(\infty))]
\right|^2
\right]
&\leq
\frac{C}{M},
&
\mathbb E\!\left[
\left|
x_{\boldsymbol a}(\boldsymbol W(\infty))
-\mathbb E[x_{\boldsymbol a}(\boldsymbol W(\infty))]
\right|^4
\right]
&\leq
\frac{C}{M^2}.
\label{eq:fa-centered-moments}
\end{align}
The bounds in \eqref{eq:fa-centered-moments} establish
\eqref{eq:fa-centered-goal}. So the second- and fourth-moment bounds in the lemma hold.
\end{proof}

\input{so-lemma-2}

\input{so-thm-1}

\input{so-thm-2}

%% file: so-lemma-2.tex
\subsection{Structural Modeling Conditions for A8-A9}
\label{subsec:structural-modeling}

A8-A9 are not easy to directly verify. Next, we will identify structural properties of the CTMC that guarantee A8-A9. These structural conditions concern only the functional form of the drift and
the instantaneous covariance-rate matrix and can therefore be verified
directly from the transition rules and the selected off-lattice extension. We verify these conditions for the power-of-two-choices model in Section \ref{sec:po2-b2}. 

\begin{itemize}

\item \textbf{S1 -- Local approximately affine drift.}
There exist constants $\eta_A>0$ and $1\leq C_A<\infty$, independent of $M$
and the realized parameter vector $\boldsymbol\theta$, and moduli
$\omega_A(s)\downarrow0$ and $\omega_r(s)\downarrow0$ as $s\downarrow0$,
such that the following holds. For each realized equilibrium interaction
input $\boldsymbol x_\chi^*$, there exist a matrix-valued function
$\boldsymbol A(\boldsymbol x;\theta)\in\mathbb R^{(n-1)\times(n-1)}$,
a vector-valued function
$\boldsymbol a(\boldsymbol x;\theta)\in\mathbb R^{n-1}$,
and a residual
$\boldsymbol r(\boldsymbol v,\boldsymbol x;\theta)\in\mathbb R^{n-1}$
such that, for every $\boldsymbol v\in{\cal D}_{n-1}$,
$\theta\in\Theta$, and $\boldsymbol x\in{\cal K}_\chi$ satisfying
$\|\boldsymbol x-\boldsymbol x_\chi^*\|\leq\eta_A$,
\begin{align}
f(\boldsymbol v,\Psi(\boldsymbol x);\theta)
&=
\boldsymbol A(\boldsymbol x;\theta)\boldsymbol v
+\boldsymbol a(\boldsymbol x;\theta)
+\boldsymbol r(\boldsymbol v,\boldsymbol x;\theta).
\label{eq:S1-local-affine-decomposition}
\end{align}
The first and second derivatives of $\boldsymbol A$ and $\boldsymbol a$
with respect to $\boldsymbol x$ satisfy
\begin{align}
\sup_{\substack{\theta\in\Theta,\ \boldsymbol x\in{\cal K}_\chi\\
\|\boldsymbol x-\boldsymbol x_\chi^*\|\leq\eta_A}}
\left(
\|D\boldsymbol A(\boldsymbol x;\theta)\|
+\|D^2\boldsymbol A(\boldsymbol x;\theta)\|
+\|D\boldsymbol a(\boldsymbol x;\theta)\|
+\|D^2\boldsymbol a(\boldsymbol x;\theta)\|
\right)
&\leq C_A,
\label{eq:S1-local-derivative-bound}
\end{align}
and their second derivatives satisfy the common modulus-of-continuity bound
\begin{align}
\sup_{\substack{\theta\in\Theta,\ \boldsymbol x,\boldsymbol y\in{\cal K}_\chi\\
\|\boldsymbol x-\boldsymbol x_\chi^*\|,
\|\boldsymbol y-\boldsymbol x_\chi^*\|\leq\eta_A\\
\|\boldsymbol x-\boldsymbol y\|\leq s}}
\left(
\|D^2\boldsymbol A(\boldsymbol x;\theta)
-D^2\boldsymbol A(\boldsymbol y;\theta)\|
+
\|D^2\boldsymbol a(\boldsymbol x;\theta)
-D^2\boldsymbol a(\boldsymbol y;\theta)\|
\right)
&\leq \omega_A(s).
\label{eq:S1-local-affine-modulus}
\end{align}
Moreover, for every $\boldsymbol x\in{\cal K}_\chi$ satisfying
$\|\boldsymbol x-\boldsymbol x_\chi^*\|\leq\eta_A$,
\begin{align}
\sup_{\boldsymbol v\in{\cal D}_{n-1},\,\theta\in\Theta}
\|\boldsymbol r(\boldsymbol v,\boldsymbol x;\theta)\|
&\leq
\omega_r(\|\boldsymbol x-\boldsymbol x_\chi^*\|)
\|\boldsymbol x-\boldsymbol x_\chi^*\|^2.
\label{eq:S1-local-affine-residual}
\end{align}
{\bf Sharper condition:} For the sharper mean-refinement rate, assume in addition that
\begin{align}
\omega_A(s)+\omega_r(s)
&\leq C_A s,
\qquad
0\leq s\leq\eta_A.
\label{eq:S1-local-affine-lipschitz}
\end{align}

\item \textbf{S2 -- Finite weighted-occupancy dependence of the noise.}
Define
\begin{align}
V_{\boldsymbol b}(\boldsymbol w)
&=
\frac1M\boldsymbol b^\top\boldsymbol Q(\boldsymbol w)\boldsymbol b
\label{eq:S2-V-A-def}
\end{align}
There exist constants $K_Q,C_Q<\infty$ and $\eta_Q>0$, independent of
$M$ and $\boldsymbol\theta$, such that for every
$\|\boldsymbol b\|_{\max}\le1$ there are block vectors
$\boldsymbol c_{1,\boldsymbol b},\ldots,\boldsymbol c_{K_Q,\boldsymbol b}$
and a twice continuously differentiable function
$\widehat V_{\boldsymbol b}:\mathbb R^{K_Q}\to\mathbb R$ satisfying
\begin{align}
V_{\boldsymbol b}(\boldsymbol w)
&=
\widehat V_{\boldsymbol b}\!\left(
\bar{x}_{\boldsymbol c_{1,\boldsymbol b}}(\boldsymbol w),\ldots,
\bar{x}_{\boldsymbol c_{K_Q,\boldsymbol b}}(\boldsymbol w)
\right)
\label{eq:S2-representation}
\end{align}
whenever
$\max_k|\bar{\boldsymbol x}_{\boldsymbol c_{k,\boldsymbol b}}(\boldsymbol w)|
\le\eta_Q$ and
\begin{align}
\max_{1\le k\le K_Q}\|\boldsymbol c_{k,\boldsymbol b}\|_{\max}
+
|\widehat V_{\boldsymbol b}(\boldsymbol0)|
+
\sup_{\|\boldsymbol z\|_\infty\le\eta_Q}
\left(
\|D\widehat V_{\boldsymbol b}(\boldsymbol z)\|
+
\|D^2\widehat V_{\boldsymbol b}(\boldsymbol z)\|
\right)
&\le C_Q.
\label{eq:S2-uniform-bound}
\end{align}
The vectors $\boldsymbol c_{k,\boldsymbol b}$ and the function
$\widehat V_{\boldsymbol b}$ may depend on $M$, $\boldsymbol\theta$, and
$\boldsymbol b$, but $K_Q$, $C_Q$, and $\eta_Q$ are uniform.

\end{itemize}

\begin{remark}
S1 does not make the full mean-field system linear. Rather, it requires the
local drift to be affine in the state of one agent when the interaction input is fixed at its equilibrium value; the full system can remain
nonlinear through $\boldsymbol x_\chi(\boldsymbol w)$. S2 is the analogous
structural condition for the noise: near equilibrium, the instantaneous
variance rate in any bounded macroscopic direction depends smoothly on only
$O(1)$ weighted occupancy deviations.     
\end{remark}

\begin{lemma}[Verification of A8--A9 from structural conditions]
Under A1--A7 and the structural modeling conditions S1--S2, A8 and A9 hold. If the sharper condition in S1 holds, then the $O(M^{-1/2})$ version of the remainder bound in A9 also holds.
\label{lem:structural-A8-A9}
\end{lemma}

\begin{proof}
Condition on $\boldsymbol\theta$ and suppress this conditioning. We verify A8 and A9 in turn. Throughout, $C$ denotes a finite generic constant, independent of $M$, $\boldsymbol\theta$, and $\boldsymbol b$, whose value may change from line to line. Write $\boldsymbol A_m(\boldsymbol x)=\boldsymbol A(\boldsymbol x;\theta_m)$,
$\boldsymbol a_m(\boldsymbol x)=\boldsymbol a(\boldsymbol x;\theta_m)$, and
$\boldsymbol r_m(\boldsymbol v,\boldsymbol x)=\boldsymbol r(\boldsymbol v,\boldsymbol x;\theta_m)$. Define the interaction-occupancy deviation
\begin{align}
\bar{\boldsymbol x}_\chi(\boldsymbol w)
&=\boldsymbol x_\chi(\boldsymbol w)-\boldsymbol x_\chi^*.
\label{eq:interaction-deviation}
\end{align} 
Lemma~\ref{lem:uniform-block-moments} gives
\begin{align}
\mathbb E\|\bar{\boldsymbol x}_\chi(\boldsymbol W(\infty))\|^2
&=O(M^{-1}),
&
\mathbb E\|\bar{\boldsymbol x}_\chi(\boldsymbol W(\infty))\|^4
&=O(M^{-2}).
\label{eq:interaction-deviation-moments}
\end{align}

Define the local S1 event
\begin{align}
{\cal B}_A
&=
\{\|\bar{\boldsymbol x}_\chi(\boldsymbol W(\infty))\|\leq\eta_A\}.
\label{eq:lemma2-BA}
\end{align}
By \eqref{eq:interaction-deviation-moments} and Markov's inequality,
\begin{align}
\mathbb P({\cal B}_A^c)
&\leq
\frac{\mathbb E\|\bar{\boldsymbol x}_\chi(\boldsymbol W(\infty))\|^4}{\eta_A^4}
=O(M^{-2}).
\label{eq:lemma2-BA-prob}
\end{align}

We first show that the residual $\boldsymbol r(\boldsymbol v,\boldsymbol x;\theta)$ does not affect the equilibrium linearization
or Hessian. Write
$\bar{\boldsymbol w}=\boldsymbol w-\boldsymbol w^*$ and
$\bar {\boldsymbol w}_m=\boldsymbol w_m-\boldsymbol w_m^*$.
By A3 and S1, $\boldsymbol r_m$ is twice continuously differentiable in
$(\boldsymbol v,\boldsymbol x)$. For every fixed
$\boldsymbol v\in{\cal D}_{n-1}$, the residual bound in S1 gives
\begin{align}
\frac{\|\boldsymbol r_m(\boldsymbol v,\boldsymbol x_\chi^*+\boldsymbol h)\|}
{\|\boldsymbol h\|^2}
&\leq \omega_r(\|\boldsymbol h\|)\longrightarrow0,
\qquad \boldsymbol h\to\boldsymbol0.
\label{eq:lemma2-residual-small}
\end{align}
Since $\|\chi\|_\infty>0$, the convex set ${\cal K}_\chi$ has nonempty
interior. Taking the limit through
$\boldsymbol x_\chi^*+\boldsymbol h\in{\cal K}_\chi$, Taylor expansion
therefore implies that $\boldsymbol r_m(\boldsymbol v,\boldsymbol x_\chi^*)$
and its first two derivatives with respect to $\boldsymbol x$ vanish. Since
$\boldsymbol r_m(\boldsymbol v,\boldsymbol x_\chi^*)=\boldsymbol0$ for every
$\boldsymbol v$, its first two $\boldsymbol v$-derivatives also vanish; and
since $D_{\boldsymbol x}\boldsymbol r_m(\boldsymbol v,\boldsymbol x_\chi^*)
=\boldsymbol0$ for every $\boldsymbol v$, the mixed second derivatives vanish
as well. Thus, with $D$ denoting differentiation with respect to the pair
$(\boldsymbol v,\boldsymbol x)$,
\begin{align}
\boldsymbol r_m(\boldsymbol v,\boldsymbol x_\chi^*)
&=\boldsymbol0,
&
D\boldsymbol r_m(\boldsymbol v,\boldsymbol x_\chi^*)
&=\boldsymbol0,
&
D^2\boldsymbol r_m(\boldsymbol v,\boldsymbol x_\chi^*)
&=\boldsymbol0.
\label{eq:lemma2-residual-jet}
\end{align}
Consequently, the residual contributes neither to the equilibrium identity,
the Jacobian $\boldsymbol J$, nor the Hessian contraction $H$. In particular,
$\boldsymbol A_m(\boldsymbol x_\chi^*)\boldsymbol w_m^*
+\boldsymbol a_m(\boldsymbol x_\chi^*)=\boldsymbol0$, and
\begin{align}
[\boldsymbol J\bar{\boldsymbol w}]_m
&=\boldsymbol A_m(\boldsymbol x_\chi^*)\bar{\boldsymbol w}_m
+D\boldsymbol A_m(\boldsymbol x_\chi^*)[\bar{\boldsymbol x}_{\chi}(\boldsymbol w)]
 \boldsymbol w_m^*
+D\boldsymbol a_m(\boldsymbol x_\chi^*)[\bar{\boldsymbol x}_{\chi}(\boldsymbol w)].
\label{eq:lemma2-J-local}
\end{align}

On ${\cal B}_A$, S1 gives
\begin{align}
\boldsymbol F_m(\boldsymbol W(\infty))
&=
\boldsymbol A_m({\boldsymbol x}_\chi(\boldsymbol W(\infty)))\boldsymbol W_m(\infty)
+\boldsymbol a_m({\boldsymbol x}_\chi(\boldsymbol W(\infty)))
+\boldsymbol r_m(\boldsymbol W_m(\infty),{\boldsymbol x}_\chi(\boldsymbol W(\infty))).
\label{eq:lemma2-local-drift}
\end{align}
More generally, for any $\boldsymbol w$ satisfying
$\|\bar{\boldsymbol x}_\chi(\boldsymbol w)\|\leq\eta_A$, define the first-order Taylor
remainders
\begin{align}
\boldsymbol R_{A,m}(\boldsymbol w)
&=
\boldsymbol A_m(\boldsymbol x_\chi^*+\bar{\boldsymbol x}_\chi(\boldsymbol w))
-\boldsymbol A_m(\boldsymbol x_\chi^*)
-D\boldsymbol A_m(\boldsymbol x_\chi^*)[\bar{\boldsymbol x}_\chi(\boldsymbol w)],
\nonumber\\
\boldsymbol R_{a,m}(\boldsymbol w)
&=
\boldsymbol a_m(\boldsymbol x_\chi^*+\bar{\boldsymbol x}_\chi(\boldsymbol w))
-\boldsymbol a_m(\boldsymbol x_\chi^*)
-D\boldsymbol a_m(\boldsymbol x_\chi^*)[\bar{\boldsymbol x}_\chi(\boldsymbol w)].
\label{eq:lemma2-first-order-remainders}
\end{align}
Taylor's theorem with integral remainder gives
\begin{align}
\boldsymbol R_{A,m}(\boldsymbol w)
&=
\int_0^1(1-s)D^2\boldsymbol A_m
\!\left(\boldsymbol x_\chi^*+s\bar{\boldsymbol x}_\chi(\boldsymbol w)\right)
[\bar{\boldsymbol x}_\chi(\boldsymbol w),\bar{\boldsymbol x}_\chi(\boldsymbol w)]\,ds,
\nonumber\\
\boldsymbol R_{a,m}(\boldsymbol w)
&=
\int_0^1(1-s)D^2\boldsymbol a_m
\!\left(\boldsymbol x_\chi^*+s\bar{\boldsymbol x}_\chi(\boldsymbol w)\right)
[\bar{\boldsymbol x}_\chi(\boldsymbol w),\bar{\boldsymbol x}_\chi(\boldsymbol w)]\,ds.
\label{eq:lemma2-first-order-integral}
\end{align}
Since the segment from $\boldsymbol x_\chi^*$ to
$\boldsymbol x_\chi^*+\bar{\boldsymbol x}_\chi(\boldsymbol w)$ stays inside the local S1
neighborhood, the uniform second-derivative bound in \eqref{eq:S1-local-derivative-bound} of S1 implies
\begin{align}
\|\boldsymbol R_{A,m}(\boldsymbol w)\|
+\|\boldsymbol R_{a,m}(\boldsymbol w)\|
&\leq C\|\bar{\boldsymbol x}_\chi(\boldsymbol w)\|^2.
\label{eq:lemma2-first-order-remainder-bound}
\end{align} for $C=C_A/2.$ Recall $\bar{\boldsymbol w}_m=\boldsymbol w_m-\boldsymbol w_m^*$. Substituting
\eqref{eq:lemma2-first-order-remainders} into the S1 decomposition and using
$\boldsymbol A_m(\boldsymbol x_\chi^*)\boldsymbol w_m^*
+\boldsymbol a_m(\boldsymbol x_\chi^*)=\boldsymbol0$ yields
\begin{align}
[\boldsymbol R_F(\boldsymbol w)]_m
&=
D\boldsymbol A_m(\boldsymbol x_\chi^*)[\bar{\boldsymbol x}_\chi(\boldsymbol w)]\bar{\boldsymbol w}_m
+\boldsymbol R_{A,m}(\boldsymbol w)\boldsymbol w_m
+\boldsymbol R_{a,m}(\boldsymbol w)
+\boldsymbol r_m(\boldsymbol w_m,\boldsymbol x_\chi^*+\bar{\boldsymbol x}_\chi(\boldsymbol w)).
\label{eq:lemma2-RF-expansion}
\end{align}
For $1\leq i\leq n-1$ and $\|\boldsymbol b\|_{\max}\leq1$, define the block
vector $\boldsymbol c^R_{i,\boldsymbol b}$ by
$\boldsymbol c^R_{i,\boldsymbol b,m}
=(D_i\boldsymbol A_m(\boldsymbol x_\chi^*))^\top\boldsymbol b_m$.
By \eqref{eq:S1-local-derivative-bound} in S1 and $\boldsymbol x_\chi^*\in{\cal K}_\chi$, 
\begin{align}
\max_{1\leq i\leq n-1}\sup_{1\leq m\leq M}
\|D_i\boldsymbol A_m(\boldsymbol x_\chi^*)\|_{\rm op}
&\leq C_A,
&
\|\boldsymbol c^R_{i,\boldsymbol b}\|_{\max}
&\leq C_A.
\label{eq:lemma2-cR-bound}
\end{align}
Define
$\overline{\boldsymbol c}^R_{i,\boldsymbol b}
=\boldsymbol c^R_{i,\boldsymbol b}/C_A$. Then
$\|\overline{\boldsymbol c}^R_{i,\boldsymbol b}\|_{\max}\leq1$ and
$\bar{x}_{\boldsymbol c^R_{i,\boldsymbol b}}(\boldsymbol W(\infty))
=C_A \bar{x}_{\overline{\boldsymbol c}^R_{i,\boldsymbol b}}(\boldsymbol W(\infty))$. Hence Lemma~\ref{lem:uniform-block-moments} gives,
uniformly in $i$ and $\boldsymbol b$,
\begin{align}
\left|\mathbb E \left[\bar{x}_{\boldsymbol c^R_{i,\boldsymbol b}}(\boldsymbol W(\infty))\right]\right|
&\leq \frac{C}{M},
&
\mathbb E|\bar{x}_{\boldsymbol c^R_{i,\boldsymbol b}}(\boldsymbol W(\infty))|^2
&\leq \frac{C}{M},
&
\mathbb E|\bar{x}_{\boldsymbol c^R_{i,\boldsymbol b}}(\boldsymbol W(\infty))|^4
&\leq \frac{C}{M^2}.
\label{eq:lemma2-cR-moments}
\end{align}
By \eqref{eq:lemma2-first-order-remainder-bound}, the residual bound in S1,
and boundedness of the single-agent state space, on ${\cal B}_A$ we have
\begin{align}
\frac1M|\boldsymbol b^\top\boldsymbol R_F(\boldsymbol W(\infty))|
&\leq
C\left(
\|\bar{\boldsymbol x}_\chi(\boldsymbol W(\infty))\|
\sum_{i=1}^{n-1}|\bar{x}_{\boldsymbol c^R_{i,\boldsymbol b}}(\boldsymbol W(\infty))|
+\left[1+\omega_r(\|\bar{\boldsymbol x}_\chi(\boldsymbol W(\infty))\|)\right]\|\bar{\boldsymbol x}_\chi(\boldsymbol W(\infty))\|^2
\right).
\label{eq:lemma2-RF-bound}
\end{align}
Since $\eta_A$ is fixed, $\omega_r$ is bounded on $[0,\eta_A]$, so the factor
$\omega_r(\|\bar{\boldsymbol x}_\chi(\boldsymbol W(\infty))\|)$ can be absorbed into $C$ on
${\cal B}_A$. Since $\boldsymbol b^\top\boldsymbol U=M \bar{x}_{\boldsymbol b}(\boldsymbol W(\infty))$,
separating the expectation over ${\cal B}_A$ and ${\cal B}_A^c$ gives
\begin{align}
\frac1M|\boldsymbol b^\top\boldsymbol E_R\boldsymbol b|
&\leq
2\left|\mathbb E\!\left[
(\boldsymbol b^\top\boldsymbol R_F)\bar{x}_{\boldsymbol b}(\boldsymbol W(\infty))
\mathbf1_{{\cal B}_A}\right]\right|
+2\mathbb E\!\left[
|\boldsymbol b^\top\boldsymbol R_F|\,|\bar{x}_{\boldsymbol b}(\boldsymbol W(\infty))|
\mathbf1_{{\cal B}_A^c}\right].
\label{eq:lemma2-ER-split}
\end{align}
For the first term, \eqref{eq:lemma2-RF-bound},
\eqref{eq:interaction-deviation-moments}, \eqref{eq:lemma2-cR-moments},
Lemma~\ref{lem:uniform-block-moments} applied to $\boldsymbol b$, and H\"older's inequality give
\begin{align}
&M\mathbb E[\|\bar{\boldsymbol x}_\chi(\boldsymbol W(\infty))\|
|\bar{x}_{\boldsymbol c^R_{i,\boldsymbol b}}(\boldsymbol W(\infty))||\bar{x}_{\boldsymbol b}(\boldsymbol W(\infty))|]\nonumber\\
\leq &
M(\mathbb E\|\bar{\boldsymbol x}_\chi(\boldsymbol W(\infty))\|^4)^{1/4}
(\mathbb E|\bar{x}_{\boldsymbol c^R_{i,\boldsymbol b}}(\boldsymbol W(\infty))|^4)^{1/4}
(\mathbb E|\bar{x}_{\boldsymbol b}(\boldsymbol W(\infty))|^2)^{1/2}
=O(M^{-1/2}),
\nonumber\\
&M\mathbb E[\|\bar{\boldsymbol x}_\chi(\boldsymbol W(\infty))\|^2|\bar{x}_{\boldsymbol b}(\boldsymbol W(\infty))|]\nonumber\\
\leq&
M(\mathbb E\|\bar{\boldsymbol x}_\chi(\boldsymbol W(\infty))\|^4)^{1/2}
(\mathbb E|\bar{x}_{\boldsymbol b}(\boldsymbol W(\infty))|^2)^{1/2}
=O(M^{-1/2}).
\label{eq:lemma2-ER-holder}
\end{align}
For the contribution on ${\cal B}_A^c$, recall that $\boldsymbol x_\chi$
is linear and $\boldsymbol F_m(\boldsymbol w^*)=\boldsymbol 0$.
Thus $[\boldsymbol R_F(\boldsymbol w)]_m$ is the first-order Taylor
remainder of $(\boldsymbol v,\boldsymbol x)\mapsto
\boldsymbol f(\boldsymbol v,\Psi(\boldsymbol x);\theta_m)$
about $(\boldsymbol w_m^*,\boldsymbol x_\chi^*)$, evaluated at
$(\boldsymbol w_m,\boldsymbol x_\chi(\boldsymbol w))$.
By A3, the second derivatives of this map are uniformly bounded on
the fixed bounded convex set ${\cal D}_{n-1}\times{\cal K}_\chi$.
Hence $\|[\boldsymbol R_F(\boldsymbol w)]_m\|\le C$ uniformly in
$m$, $M$, $\boldsymbol\theta$, and $\boldsymbol w\in{\cal D}$.

For $\|\boldsymbol b\|_{\max}\le1$, summing over $m$ gives
$|\boldsymbol b^\top\boldsymbol R_F(\boldsymbol w)|\le CM$.
The simplex constraint also gives $|\bar x_{\boldsymbol b}(\boldsymbol w)|
\le M^{-1}\sum_{m=1}^M\|\boldsymbol w_m-\boldsymbol w_m^*\|\le2$.
Since $\mathbb P({\cal B}_A^c)=O(M^{-2})$, we obtain
\begin{align}
\mathbb E\!\left[
  \bigl|\boldsymbol b^\top\boldsymbol R_F(\boldsymbol W(\infty))\bigr|\,
  \bigl|\bar x_{\boldsymbol b}(\boldsymbol W(\infty))\bigr|\,
  \mathbf1_{{\cal B}_A^c}
\right]
&\le 2CM\,\mathbb P({\cal B}_A^c)=O(M^{-1}).
\label{eq:lemma2-ER-bad-event}
\end{align}
Combining \eqref{eq:lemma2-ER-split}--\eqref{eq:lemma2-ER-bad-event}
verifies the $\boldsymbol E_R$ bound in A8.

We next verify the $\boldsymbol E_Q$ bound in A8. For $\|\boldsymbol b\|_{\max}\leq1$, every feasible
jump changes at most $c_e$ blocks, and the distance between two reduced one-hot
states in one block is at most $\sqrt2$. Hence
$|\boldsymbol b^\top\boldsymbol\Delta_r|\leq\sqrt2\,c_e$, and A1 gives
\begin{align}
0\leq V_{\boldsymbol b}(\boldsymbol W)
&=\frac1M\sum_{r\in{\cal R}}q_r^{\rm CTMC}(\boldsymbol W)
(\boldsymbol b^\top\boldsymbol\Delta_r)^2
\leq 2c c_e^2,
\qquad \boldsymbol W\in{\cal W}.
\label{eq:lemma2-V-discrete-bound}
\end{align}
At equilibrium, S2 gives
\begin{align}
V_{\boldsymbol b}(\boldsymbol w^*)
&=\widehat V_{\boldsymbol b}(\boldsymbol 0),
&
\sup_{\|\boldsymbol b\|_{\max}\leq1}|V_{\boldsymbol b}(\boldsymbol w^*)|
&\leq C_Q.
\label{eq:lemma2-V-equilibrium-bound}
\end{align}
For fixed $\boldsymbol b$, let
$\boldsymbol c_{1,\boldsymbol b},\ldots,\boldsymbol c_{K_Q,\boldsymbol b}$ and
$\widehat V_{\boldsymbol b}$ be given by S2, and define event 
\begin{align}
{\cal B}_{\boldsymbol b}
&=\left\{\max_{1\leq k\leq K_Q}
|\bar{\boldsymbol x}_{\boldsymbol c_{k,\boldsymbol b}}({\boldsymbol W}(\infty))|\leq\eta_Q\right\}.
\label{eq:lemma2-Bb}
\end{align} We will use shorthand notation $\bar{\boldsymbol x}_{\boldsymbol c_{k,\boldsymbol b}}=\bar{\boldsymbol x}_{\boldsymbol c_{k,\boldsymbol b}}({\boldsymbol W}(\infty))$ in the rest of the proof.
Separating the expectation over ${\cal B}_{\boldsymbol b}$ and
${\cal B}_{\boldsymbol b}^c$,  we obtain
\begin{align}
\frac1M|\boldsymbol b^\top\boldsymbol E_Q\boldsymbol b|
&=\left|\mathbb E\!\left[V_{\boldsymbol b}(\boldsymbol W(\infty))
-V_{\boldsymbol b}(\boldsymbol w^*)\right]\right|
\nonumber\\
&\leq
\left|\mathbb E\!\left[
\bigl(V_{\boldsymbol b}(\boldsymbol W(\infty))-V_{\boldsymbol b}(\boldsymbol w^*)\bigr)
\mathbf 1_{{\cal B}_{\boldsymbol b}}\right]\right|
+\mathbb E\!\left[
\left|V_{\boldsymbol b}(\boldsymbol W(\infty))-V_{\boldsymbol b}(\boldsymbol w^*)\right|
\mathbf 1_{{\cal B}_{\boldsymbol b}^c}\right].
\nonumber
\end{align}

On ${\cal B}_{\boldsymbol b}$, define
$\boldsymbol z_{\boldsymbol b}
=(\bar{\boldsymbol x}_{\boldsymbol c_{1,\boldsymbol b}},\ldots,
\bar{\boldsymbol x}_{\boldsymbol c_{K_Q,\boldsymbol b}})^\top$.
Then $\|\boldsymbol z_{\boldsymbol b}\|_\infty\leq\eta_Q$, so the S2 representation is valid at
$\boldsymbol W(\infty)$. Moreover,
$\|s\boldsymbol z_{\boldsymbol b}\|_\infty\leq\eta_Q$ for every $s\in[0,1]$,
so the S2 bound on $D^2\widehat V_{\boldsymbol b}$ holds along the entire segment
from $\boldsymbol 0$ to $\boldsymbol z_{\boldsymbol b}$. Since
$\bar{\boldsymbol x}_{\boldsymbol c_{k,\boldsymbol b}}(\boldsymbol w^*)=0$ for every $k$,
$V_{\boldsymbol b}(\boldsymbol w^*)=\widehat V_{\boldsymbol b}(\boldsymbol 0)$.
Taylor's theorem at $\boldsymbol 0$ therefore gives
\begin{align}
V_{\boldsymbol b}(\boldsymbol W(\infty))-V_{\boldsymbol b}(\boldsymbol w^*)
&=
\widehat V_{\boldsymbol b}(\boldsymbol z_{\boldsymbol b})
-\widehat V_{\boldsymbol b}(\boldsymbol 0)
=\sum_{k=1}^{K_Q}\kappa_{k,\boldsymbol b}
\bar{\boldsymbol x}_{\boldsymbol c_{k,\boldsymbol b}}+R_{\boldsymbol b},
\end{align}
where
\begin{align}
\kappa_{k,\boldsymbol b}
=\frac{\partial \widehat V_{\boldsymbol b}}{\partial z_k}(\boldsymbol 0)\qquad \hbox{and}\qquad 
R_{\boldsymbol b}
&=\int_0^1(1-s)\boldsymbol z_{\boldsymbol b}^\top
D^2\widehat V_{\boldsymbol b}(s\boldsymbol z_{\boldsymbol b})
\boldsymbol z_{\boldsymbol b}\,ds.
\label{eq:lemma2-V-taylor}
\end{align}

On ${\cal B}_{\boldsymbol b}^c$, \eqref{eq:lemma2-V-discrete-bound}--
\eqref{eq:lemma2-V-equilibrium-bound} give a uniform bound on
$|V_{\boldsymbol b}(\boldsymbol W(\infty))-V_{\boldsymbol b}(\boldsymbol w^*)|$.
By Lemma~\ref{lem:uniform-block-moments}, the uniform block-max bound in S2, a union bound, and Markov's inequality applied to
$|\bar{\boldsymbol x}_{\boldsymbol c_{k,\boldsymbol b}}|^4$,
\begin{align}
\mathbb P({\cal B}_{\boldsymbol b}^c)
&\leq
\sum_{k=1}^{K_Q}
\mathbb P\left(|\bar{\boldsymbol x}_{\boldsymbol c_{k,\boldsymbol b}}|>\eta_Q\right)
\leq
\sum_{k=1}^{K_Q}
\frac{\mathbb E|\bar{\boldsymbol x}_{\boldsymbol c_{k,\boldsymbol b}}|^4}{\eta_Q^4}
=O(M^{-2}).
\label{eq:lemma2-Bb-prob}
\end{align}
Therefore, we have
\begin{align}
\frac1M|\boldsymbol b^\top\boldsymbol E_Q\boldsymbol b|
&\leq
\sum_{k=1}^{K_Q}|\kappa_{k,\boldsymbol b}|
\left|\mathbb E\!\left[\bar{\boldsymbol x}_{\boldsymbol c_{k,\boldsymbol b}}
\mathbf 1_{{\cal B}_{\boldsymbol b}}\right]\right|
+\mathbb E\!\left[|R_{\boldsymbol b}|\mathbf 1_{{\cal B}_{\boldsymbol b}}\right]
+O(M^{-2}).
\end{align}

Next, by S2,
\begin{align}
|\kappa_{k,\boldsymbol b}|
&=
\left|\frac{\partial \widehat V_{\boldsymbol b}}{\partial z_k}(\boldsymbol 0)\right|
\leq
\|D\widehat V_{\boldsymbol b}(\boldsymbol 0)\|
\leq C_Q,
\qquad 1\leq k\leq K_Q.
\label{eq:lemma2-kappa-bound}
\end{align}
Moreover, by Lemma~\ref{lem:uniform-block-moments} and \eqref{eq:lemma2-Bb-prob},
\begin{align}
\left|\mathbb E[\bar{\boldsymbol x}_{\boldsymbol c_{k,\boldsymbol b}}\mathbf 1_{{\cal B}_{\boldsymbol b}}]\right|
&\leq |\mathbb E \bar{\boldsymbol x}_{\boldsymbol c_{k,\boldsymbol b}}|
+(\mathbb E|\bar{\boldsymbol x}_{\boldsymbol c_{k,\boldsymbol b}}|^2)^{1/2}
\mathbb P({\cal B}_{\boldsymbol b}^c)^{1/2}
=O(M^{-1}).
\label{eq:lemma2-linear-on-B}
\end{align}
From the definition of $R_{\boldsymbol b},$ we have 
\begin{align}
|R_{\boldsymbol b}|
\leq \frac{C_Q}{2}\sum_{k=1}^{K_Q}
|\bar{\boldsymbol x}_{\boldsymbol c_{k,\boldsymbol b}}|^2.
\end{align}
By \eqref{eq:lemma2-V-taylor} and Lemma~\ref{lem:uniform-block-moments},
\begin{align}
\mathbb E[|R_{\boldsymbol b}|\mathbf 1_{{\cal B}_{\boldsymbol b}}]
&\leq \frac{C_Q}{2}\sum_{k=1}^{K_Q}
\mathbb E\!
\left[|\bar{\boldsymbol x}_{\boldsymbol c_{k,\boldsymbol b}}|^2
\mathbf 1_{{\cal B}_{\boldsymbol b}}\right]
\leq \frac{C_Q}{2}\sum_{k=1}^{K_Q}
\mathbb E|\bar{\boldsymbol x}_{\boldsymbol c_{k,\boldsymbol b}}|^2
=O(M^{-1}).
\label{eq:lemma2-remainder-on-B}
\end{align}

Combining the results above, we conclude that 
\begin{align}
\frac1M|\boldsymbol b^\top\boldsymbol E_Q\boldsymbol b|=O(M^{-1}).
\label{eq:lemma2-EQ}
\end{align}
uniformly over $\|\boldsymbol b\|_{\max}\leq1$. This is in fact the stronger $O(M^{-1})$ estimate for $\boldsymbol E_Q$, and hence verifies the $O(M^{-1/2})$ bound required in A8. Thus A8 holds.

We next verify the Hessian-control requirement in A9 using two full-state directions
$\boldsymbol u,\boldsymbol v\in\mathbb R^N$. By \eqref{eq:lemma2-residual-jet}, the residual makes no contribution to the
Hessian contraction $H$. Since $\boldsymbol x_\chi$ is linear,
for every full-state direction $\boldsymbol u$,
\begin{align}
D\boldsymbol x_\chi[\boldsymbol u]
&=\frac1M\sum_{m=1}^M\chi(\theta_m)\boldsymbol u_m.
\label{eq:lemma2-xchi-direction}
\end{align}
Differentiating the local S1 decomposition twice at equilibrium gives
\begin{align}
D^2\boldsymbol F_m(\boldsymbol w^*)[\boldsymbol u,\boldsymbol v]
&=D\boldsymbol A_m(\boldsymbol x_\chi^*)[D\boldsymbol x_\chi[\boldsymbol u]]\boldsymbol v_m
+D\boldsymbol A_m(\boldsymbol x_\chi^*)[D\boldsymbol x_\chi[\boldsymbol v]]\boldsymbol u_m
\nonumber\\
&\quad+D^2\boldsymbol A_m(\boldsymbol x_\chi^*)
[D\boldsymbol x_\chi[\boldsymbol u],D\boldsymbol x_\chi[\boldsymbol v]]\boldsymbol w_m^*
+D^2\boldsymbol a_m(\boldsymbol x_\chi^*)
[D\boldsymbol x_\chi[\boldsymbol u],D\boldsymbol x_\chi[\boldsymbol v]].
\label{eq:lemma2-H-bilinear}
\end{align}
For $1\leq i,j\leq n-1$, define $\boldsymbol d^{(i)}$ to be a block vector such that its $m$th block is $\boldsymbol d_m^{(i)}=\chi(\theta_m){\boldsymbol e}_i$ and define 
\begin{align}
\boldsymbol r_{i,\boldsymbol b,m}
&=(D_i\boldsymbol A_m(\boldsymbol x_\chi^*))^\top\boldsymbol b_m,
&
\rho_{ij,\boldsymbol b}
&=\frac1M\sum_{m=1}^M\boldsymbol b_m^\top
\left(D^2_{ij}\boldsymbol A_m(\boldsymbol x_\chi^*)\boldsymbol w_m^*
+D^2_{ij}\boldsymbol a_m(\boldsymbol x_\chi^*)\right).
\label{eq:lemma2-r-rho-H}
\end{align}
Recall that $\boldsymbol L_{\boldsymbol a}=M^{-1/2}\boldsymbol a^\top$. Hence
\begin{align}
D x_{\chi,i}[\boldsymbol u]
&=M^{-1}(\boldsymbol d^{(i)})^\top\boldsymbol u
=M^{-1/2}\boldsymbol L_{\boldsymbol d^{(i)}}\boldsymbol u,
&
(\boldsymbol r_{i,\boldsymbol b})^\top\boldsymbol u
&=\sqrt M\,\boldsymbol L_{\boldsymbol r_{i,\boldsymbol b}}\boldsymbol u.
\label{eq:lemma2-L-relations}
\end{align}
Thus, summing \eqref{eq:lemma2-H-bilinear} against $\boldsymbol b$ gives
\begin{align}
\boldsymbol b^\top D^2\boldsymbol F(\boldsymbol w^*)[\boldsymbol u,\boldsymbol v]
&=\sum_{i=1}^{n-1}\left[
(\boldsymbol L_{\boldsymbol d^{(i)}}\boldsymbol u)
(\boldsymbol L_{\boldsymbol r_{i,\boldsymbol b}}\boldsymbol v)
+(\boldsymbol L_{\boldsymbol d^{(i)}}\boldsymbol v)
(\boldsymbol L_{\boldsymbol r_{i,\boldsymbol b}}\boldsymbol u)
\right]
+\sum_{i,j=1}^{n-1}\rho_{ij,\boldsymbol b}
(\boldsymbol L_{\boldsymbol d^{(i)}}\boldsymbol u)
(\boldsymbol L_{\boldsymbol d^{(j)}}\boldsymbol v).
\label{eq:lemma2-H-uv}
\end{align}
The symmetry of the second derivatives implies
$\rho_{ij,\boldsymbol b}=\rho_{ji,\boldsymbol b}$. Now set
\begin{align}
\boldsymbol B
&=\frac12\left(\boldsymbol u\boldsymbol v^\top
+\boldsymbol v\boldsymbol u^\top\right).
\label{eq:lemma2-B-uv}
\end{align}
By the componentwise definition of the Hessian contraction $H$, symmetry of the component Hessians gives
$H(\boldsymbol B)=D^2\boldsymbol F(\boldsymbol w^*)[\boldsymbol u,\boldsymbol v]$.
Thus \eqref{eq:lemma2-H-uv} becomes
\begin{align}
\boldsymbol b^\top H(\boldsymbol B)
&=2\sum_{i=1}^{n-1}
\boldsymbol L_{\boldsymbol r_{i,\boldsymbol b}}\boldsymbol B
\boldsymbol L_{\boldsymbol d^{(i)}}^\top
+\sum_{i,j=1}^{n-1}\rho_{ij,\boldsymbol b}
\boldsymbol L_{\boldsymbol d^{(i)}}\boldsymbol B
\boldsymbol L_{\boldsymbol d^{(j)}}^\top.
\label{eq:lemma2-H-explicit}
\end{align}
Every symmetric matrix is a linear combination of matrices of the form
$\boldsymbol u\boldsymbol u^\top$. Since both sides of
\eqref{eq:lemma2-H-explicit} are linear in $\boldsymbol B$, the identity therefore
holds for every symmetric $\boldsymbol B$.

We now put \eqref{eq:lemma2-H-explicit} in the form required by A9. By S1,
$\|\boldsymbol b\|_{\max}\leq1$, boundedness of $\chi$, and
$\boldsymbol w_m^*\in{\cal D}_{n-1}$, there exists $C_0<\infty$, independent of
$M$, $\boldsymbol\theta$, and $\boldsymbol b$, such that
\begin{align}
\max_{1\leq i\leq n-1}
\left\{\|\boldsymbol r_{i,\boldsymbol b}\|_{\max},
\|\boldsymbol d^{(i)}\|_{\max}\right\}
+\max_{1\leq i,j\leq n-1}|\rho_{ij,\boldsymbol b}|
&\leq C_0.
\label{eq:lemma2-H-coefficient-bound}
\end{align}
Specifically, the first two bounds follow from the uniform bound on
$D\boldsymbol A$ in S1 and $\|\boldsymbol d^{(i)}\|_{\max}\leq\|\chi\|_\infty$,
while
\begin{align}
|\rho_{ij,\boldsymbol b}|
&\leq
\frac1M\sum_{m=1}^M\|\boldsymbol b_m\|
\left(
\|D^2_{ij}\boldsymbol A_m(\boldsymbol x_\chi^*)\|\,\|\boldsymbol w_m^*\|
+\|D^2_{ij}\boldsymbol a_m(\boldsymbol x_\chi^*)\|
\right)
\leq C_0
\label{eq:lemma2-rho-bound}
\end{align}
by the uniform second-derivative bounds in S1.

There are $(n-1)$ terms in the first sum of
\eqref{eq:lemma2-H-explicit} and $(n-1)^2$ terms in the second. Set
$K_H=(n-1)+(n-1)^2=n(n-1)$ and relabel the corresponding pairs of block
vectors as
\begin{align}
(\boldsymbol r_{k,\boldsymbol b},\boldsymbol s_{k,\boldsymbol b})
&\in
\left\{
(\boldsymbol r_{i,\boldsymbol b},\boldsymbol d^{(i)}):1\leq i\leq n-1
\right\}
\cup
\left\{
(\boldsymbol d^{(i)},\boldsymbol d^{(j)}):1\leq i,j\leq n-1
\right\}.
\label{eq:lemma2-H-pairs}
\end{align}
Then \eqref{eq:lemma2-H-explicit} and the triangle inequality imply
\begin{align}
|\boldsymbol b^\top H(\boldsymbol B)|
&\leq
2\sum_{i=1}^{n-1}
\left|
\boldsymbol L_{\boldsymbol r_{i,\boldsymbol b}}\boldsymbol B
\boldsymbol L_{\boldsymbol d^{(i)}}^\top
\right|
+C_0\sum_{i,j=1}^{n-1}
\left|
\boldsymbol L_{\boldsymbol d^{(i)}}\boldsymbol B
\boldsymbol L_{\boldsymbol d^{(j)}}^\top
\right|
\nonumber\\
&\leq
C_H\sum_{k=1}^{K_H}
\left|
\boldsymbol L_{\boldsymbol r_{k,\boldsymbol b}}\boldsymbol B
\boldsymbol L_{\boldsymbol s_{k,\boldsymbol b}}^\top
\right|,
\qquad
\max_{1\leq k\leq K_H}
\left\{
\|\boldsymbol r_{k,\boldsymbol b}\|_{\max},
\|\boldsymbol s_{k,\boldsymbol b}\|_{\max}
\right\}
\leq C_H,
\label{eq:lemma2-H-final-bound}
\end{align}
for a constant $C_H<\infty$ independent of $M$, $\boldsymbol\theta$, and
$\boldsymbol b$. This verifies the Hessian-control requirement in A9.

It remains to verify the post-quadratic drift-remainder requirement in A9. The following
argument is used only on ${\cal B}_A$, where the local S1 representation is valid.
For $\boldsymbol w\in{\cal D}$ satisfying
$\|\overline{\boldsymbol x}_\chi(\boldsymbol w)\|\le\eta_A$,
Taylor expansion of $\boldsymbol A_m$ and $\boldsymbol a_m$ about
$\boldsymbol x_\chi^*$ gives
\begin{align}
\boldsymbol A_m(\boldsymbol x_\chi(\boldsymbol w))
&= \boldsymbol A_m(\boldsymbol x_\chi^*)
 + D\boldsymbol A_m(\boldsymbol x_\chi^*)[\overline{\boldsymbol x}_\chi(\boldsymbol w)]
 + \underbrace{\tfrac12 D^2\boldsymbol A_m(\boldsymbol x_\chi^*)
   [\overline{\boldsymbol x}_\chi(\boldsymbol w),\overline{\boldsymbol x}_\chi(\boldsymbol w)]
   + \widetilde{\boldsymbol R}_{A,m}(\boldsymbol w)}_{\boldsymbol R_{A,m}(\boldsymbol w)},
 \notag\\
\boldsymbol a_m(\boldsymbol x_\chi(\boldsymbol w))
&= \boldsymbol a_m(\boldsymbol x_\chi^*)
 + D\boldsymbol a_m(\boldsymbol x_\chi^*)[\overline{\boldsymbol x}_\chi(\boldsymbol w)]
 + \underbrace{\tfrac12 D^2\boldsymbol a_m(\boldsymbol x_\chi^*)
   [\overline{\boldsymbol x}_\chi(\boldsymbol w),\overline{\boldsymbol x}_\chi(\boldsymbol w)]
   + \widetilde{\boldsymbol R}_{a,m}(\boldsymbol w)}_{\boldsymbol R_{a,m}(\boldsymbol w)}.
 \label{eq:lemma2-coefficient-taylor}
\end{align}

For $\|\overline{\boldsymbol x}_\chi(\boldsymbol w)\|\le\eta_A$,
apply Taylor's integral remainder formula to $\boldsymbol A_m$ and
$\boldsymbol a_m$ along the segment from $\boldsymbol x_\chi^*$ to
$\boldsymbol x_\chi(\boldsymbol w)$. The common modulus bound in S1
and $\int_0^1(1-s)\,ds=1/2$ give
\begin{align}
\|\widetilde{\boldsymbol R}_{A,m}(\boldsymbol w)\|
+\|\widetilde{\boldsymbol R}_{a,m}(\boldsymbol w)\|
&\le \tfrac12\omega_A\!\left(\|\overline{\boldsymbol x}_\chi(\boldsymbol w)\|\right)
     \|\overline{\boldsymbol x}_\chi(\boldsymbol w)\|^2.
\label{eq:lemma2-coefficient-remainder}
\end{align}
Let $\overline{\boldsymbol w}=\boldsymbol w-\boldsymbol w^*$.
Taking the derivative
at $\boldsymbol w^*$ in the direction
$\overline{\boldsymbol w}=\boldsymbol w-\boldsymbol w^*$ gives
\begin{align}
D\boldsymbol x_\chi[\overline{\boldsymbol w}]
&=\left.\frac{d}{dt}\boldsymbol x_\chi
    (\boldsymbol w^*+t\overline{\boldsymbol w})\right|_{t=0}
 =\frac{1}{M}\sum_{m=1}^M\chi(\theta_m)
    (\boldsymbol w_m-\boldsymbol w_m^*) \notag\\
&=\boldsymbol x_\chi(\boldsymbol w)-\boldsymbol x_\chi(\boldsymbol w^*)
 =\overline{\boldsymbol x}_\chi(\boldsymbol w).
\label{eq:interaction-occupancy-centered-direction}
\end{align}
Using $D\boldsymbol x_\chi[\overline{\boldsymbol w}]
=\overline{\boldsymbol x}_\chi(\boldsymbol w)$ in
$D^2\boldsymbol F_m(\boldsymbol w^*)[\overline{\boldsymbol w},\overline{\boldsymbol w}]$ gives
\begin{align}
\tfrac12\bigl[H(\overline{\boldsymbol w}\,\overline{\boldsymbol w}^{\top})\bigr]_m
&= D\boldsymbol A_m(\boldsymbol x_\chi^*)[\overline{\boldsymbol x}_\chi(\boldsymbol w)]
   \overline{\boldsymbol w}_m
 + \tfrac12 D^2\boldsymbol A_m(\boldsymbol x_\chi^*)
   [\overline{\boldsymbol x}_\chi(\boldsymbol w),\overline{\boldsymbol x}_\chi(\boldsymbol w)]
   \boldsymbol w_m^* + \tfrac12 D^2\boldsymbol a_m(\boldsymbol x_\chi^*)
   [\overline{\boldsymbol x}_\chi(\boldsymbol w),\overline{\boldsymbol x}_\chi(\boldsymbol w)].
\label{eq:lemma2-quadratic-drift-block}
\end{align}
In the block formula for $\boldsymbol R_F(\boldsymbol w)$, expand
$\boldsymbol R_{A,m}(\boldsymbol w)$ and $\boldsymbol R_{a,m}(\boldsymbol w)$
using \eqref{eq:lemma2-coefficient-taylor}. Subtracting
\eqref{eq:lemma2-quadratic-drift-block} cancels the $D\boldsymbol A_m$ and
$D^2\boldsymbol a_m$ terms and leaves
$\boldsymbol w_m-\boldsymbol w_m^*=\overline{\boldsymbol w}_m$
in the $D^2\boldsymbol A_m$ term. By the definition of $\boldsymbol T_F$,
we obtain
\begin{align}
[\boldsymbol T_F(\boldsymbol w)]_m
&= \tfrac12 D^2\boldsymbol A_m(\boldsymbol x_\chi^*)
   [\overline{\boldsymbol x}_\chi(\boldsymbol w),\overline{\boldsymbol x}_\chi(\boldsymbol w)]
   \overline{\boldsymbol w}_m
 + \widetilde{\boldsymbol R}_{A,m}(\boldsymbol w)\boldsymbol w_m
 + \widetilde{\boldsymbol R}_{a,m}(\boldsymbol w)
 + \boldsymbol r_m(\boldsymbol w_m,\boldsymbol x_\chi(\boldsymbol w)).
\label{eq:lemma2-post-quadratic-drift-block}
\end{align}

Fix a block vector $\boldsymbol b$ with $\|\boldsymbol b\|_{\max}\le1$.
Let $\bar x_{\chi,i}(\boldsymbol w)$ denote the $i$th coordinate of
$\overline{\boldsymbol x}_\chi(\boldsymbol w)$. For $1\le i,j\le n-1$,
define the block vector $\widetilde{\boldsymbol c}_{ij,\boldsymbol b}$ by
$\widetilde{\boldsymbol c}_{ij,\boldsymbol b,m}
=(D^2_{ij}\boldsymbol A_m(\boldsymbol x_\chi^*))^\top\boldsymbol b_m$.
The second-derivative bound in S1 gives
$\|\widetilde{\boldsymbol c}_{ij,\boldsymbol b}\|_{\max}\le C_A$.
Expanding the Hessian in coordinates and using
$\bar x_{\boldsymbol c}(\boldsymbol w)
=M^{-1}\sum_{m=1}^M\boldsymbol c_m^\top\overline{\boldsymbol w}_m$ yields
\begin{align}
\frac{1}{2M}\sum_{m=1}^M \boldsymbol b_m^\top
 D^2\boldsymbol A_m(\boldsymbol x_\chi^*)
 [\overline{\boldsymbol x}_\chi(\boldsymbol w),
  \overline{\boldsymbol x}_\chi(\boldsymbol w)]\overline{\boldsymbol w}_m
&=\frac12\sum_{i,j=1}^{n-1}
 \bar x_{\chi,i}(\boldsymbol w)\bar x_{\chi,j}(\boldsymbol w)\,
 \bar x_{\widetilde{\boldsymbol c}_{ij,\boldsymbol b}}(\boldsymbol w).
\label{eq:lemma2-tf-coordinate}
\end{align}
Set $\omega(s)=\omega_A(s)+\omega_r(s)$. By
\eqref{eq:lemma2-coefficient-remainder}, the residual bound in S1, and
$\|\boldsymbol w_m\|\le1$, the sum of the norms of the last three terms
in \eqref{eq:lemma2-post-quadratic-drift-block} is at most
$\omega(\|\overline{\boldsymbol x}_\chi(\boldsymbol w)\|)
 \|\overline{\boldsymbol x}_\chi(\boldsymbol w)\|^2$.
Using $|\bar x_{\chi,i}(\boldsymbol w)\bar x_{\chi,j}(\boldsymbol w)|
\le\|\overline{\boldsymbol x}_\chi(\boldsymbol w)\|^2$
in \eqref{eq:lemma2-tf-coordinate}, we obtain, for
$\|\overline{\boldsymbol x}_\chi(\boldsymbol w)\|\le\eta_A$,
\begin{align}
\frac1M|\boldsymbol b^\top\boldsymbol T_F(\boldsymbol w)|
&\le\|\overline{\boldsymbol x}_\chi(\boldsymbol w)\|^2
 \left(\frac12\sum_{i,j=1}^{n-1}
 |\bar x_{\widetilde{\boldsymbol c}_{ij,\boldsymbol b}}(\boldsymbol w)|
 +\omega(\|\overline{\boldsymbol x}_\chi(\boldsymbol w)\|)\right).
\label{eq:lemma2-tf-local}
\end{align}

On ${\cal B}_A$, evaluate \eqref{eq:lemma2-tf-local} at
$\boldsymbol w=\boldsymbol W(\infty)$. Conditional on $\boldsymbol\theta$,
the coefficient vector $\widetilde{\boldsymbol c}_{ij,\boldsymbol b}$ is
deterministic. Applying Lemma~1 to
$\widetilde{\boldsymbol c}_{ij,\boldsymbol b}/C_A$ gives
$\mathbb E|\bar x_{\widetilde{\boldsymbol c}_{ij,\boldsymbol b}}
(\boldsymbol W(\infty))|^2\le C/M$ uniformly in $i,j,\boldsymbol b$.
By Cauchy--Schwarz and the fourth-moment bound for
$\overline{\boldsymbol x}_\chi(\boldsymbol W(\infty))$,
\begin{align}
&M\mathbb E\!\left[
 \|\overline{\boldsymbol x}_\chi(\boldsymbol W(\infty))\|^2
 \bigl|\bar x_{\widetilde{\boldsymbol c}_{ij,\boldsymbol b}}
       (\boldsymbol W(\infty))\bigr|
 \mathbf1_{{\cal B}_A}\right] \notag\\
&\quad\le M
 \bigl(\mathbb E\|\overline{\boldsymbol x}_\chi(\boldsymbol W(\infty))\|^4
 \bigr)^{1/2}
 \bigl(\mathbb E|\bar x_{\widetilde{\boldsymbol c}_{ij,\boldsymbol b}}
       (\boldsymbol W(\infty))|^2\bigr)^{1/2}
 =O(M^{-1/2}).
\label{eq:lemma2-tf-mixed-moment}
\end{align}
There are only $(n-1)^2$ terms, so their sum has the same order.

For the modulus term, fix $0<\eta<\eta_A$ and split ${\cal B}_A$ according
to whether $\|\overline{\boldsymbol x}_\chi(\boldsymbol W(\infty))\|\le\eta$.
When this norm exceeds $\eta$, we have
$\|\overline{\boldsymbol x}_\chi(\boldsymbol W(\infty))\|^2
\le\eta^{-2}\|\overline{\boldsymbol x}_\chi(\boldsymbol W(\infty))\|^4$.
The second- and fourth-moment bounds therefore give
\begin{align}
&M\mathbb E\!\left[
 \omega\!\left(\|\overline{\boldsymbol x}_\chi(\boldsymbol W(\infty))\|\right)
 \|\overline{\boldsymbol x}_\chi(\boldsymbol W(\infty))\|^2
 \mathbf1_{{\cal B}_A}\right] \notag\\
&\quad\le M\omega(\eta)
 \mathbb E\|\overline{\boldsymbol x}_\chi(\boldsymbol W(\infty))\|^2
 +\frac{M\omega(\eta_A)}{\eta^2}
 \mathbb E\|\overline{\boldsymbol x}_\chi(\boldsymbol W(\infty))\|^4
 \le C\omega(\eta)+\frac{C}{M\eta^2}.
\label{eq:lemma2-tf-modulus}
\end{align}
Letting $M\to\infty$ and then $\eta\downarrow0$ shows that the left-hand
side of \eqref{eq:lemma2-tf-modulus} is $o(1)$.
Under the sharper condition in S1, $\omega(s)\le C_A s$ for $0\le s\le\eta_A$,
so the fourth-moment bound gives
\begin{align}
&M\mathbb E\!\left[
 \omega\!\left(\|\overline{\boldsymbol x}_\chi(\boldsymbol W(\infty))\|\right)
 \|\overline{\boldsymbol x}_\chi(\boldsymbol W(\infty))\|^2
 \mathbf1_{{\cal B}_A}\right] \notag\\
&\quad\le C_A M
 \mathbb E\|\overline{\boldsymbol x}_\chi(\boldsymbol W(\infty))\|^3
 \le C_A M
 \bigl(\mathbb E\|\overline{\boldsymbol x}_\chi(\boldsymbol W(\infty))\|^4
 \bigr)^{3/4}
 =O(M^{-1/2}).
\label{eq:lemma2-tf-modulus-sharp}
\end{align}

On ${\cal B}_A^c$, use the bound $\|[\boldsymbol R_F(\boldsymbol w)]_m\|\le C$
on ${\cal D}$ established in the verification of A8. All derivatives in
\eqref{eq:lemma2-quadratic-drift-block} are evaluated at $\boldsymbol x_\chi^*$,
where S1 bounds them uniformly. Since $\|\overline{\boldsymbol w}_m\|\le2$
and $\|\overline{\boldsymbol x}_\chi(\boldsymbol w)\|\le2\|\chi\|_\infty$,
the quadratic term is also uniformly bounded on ${\cal D}$.
Thus $\|[\boldsymbol T_F(\boldsymbol w)]_m\|\le C$, and summing over $m$
gives $|\boldsymbol b^\top\boldsymbol T_F(\boldsymbol w)|\le CM$.
Since $\mathbb P({\cal B}_A^c)=O(M^{-2})$, we have
\begin{align}
\mathbb E\!\left[
 |\boldsymbol b^\top\boldsymbol T_F(\boldsymbol W(\infty))|
 \mathbf1_{{\cal B}_A^c}\right]
&\le CM\,\mathbb P({\cal B}_A^c)=O(M^{-1}).
\label{eq:lemma2-tf-tail}
\end{align}

The estimates \eqref{eq:lemma2-tf-local}--\eqref{eq:lemma2-tf-tail}
are uniform in $\boldsymbol b$. Taking the supremum after splitting
the expectation over ${\cal B}_A$ and ${\cal B}_A^c$ gives
\begin{align}
\sup_{\|\boldsymbol b\|_{\max}\le1}
 \left|\boldsymbol b^\top\mathbb E\!\left[
 \boldsymbol T_F(\boldsymbol W(\infty))\mid\boldsymbol\theta\right]\right|
&=o(1).
\label{eq:lemma2-tf-conclusion}
\end{align}
Under the sharper condition in S1, \eqref{eq:lemma2-tf-mixed-moment},
\eqref{eq:lemma2-tf-modulus-sharp}, and \eqref{eq:lemma2-tf-tail} improve
the bound in \eqref{eq:lemma2-tf-conclusion} to $O(M^{-1/2})$.
This proves the remaining requirement of A9 and completes the proof.
\end{proof}

%% file: so-thm-1.tex
\subsection{Covariance Approximation and Mean Refinement}

Recall that
\begin{align*}
\boldsymbol\Gamma_\varphi
=L_\varphi\overline{\boldsymbol\Sigma}L_\varphi^\top\quad\hbox{and}\quad
\boldsymbol\Sigma
=\mathbb E[
(\boldsymbol W(\infty)-\boldsymbol w^*)(\boldsymbol W(\infty)-\boldsymbol w^*)^\top
\mid\boldsymbol\theta].
\end{align*} We have the following theorem.

\begin{theorem}[Second-order covariance approximation]
\label{thm:covariance-approximation}
Assuming A1--A8, we have \begin{align}
\sup_{\substack{
\|\boldsymbol a\|_{\max}\le1\\
\|\boldsymbol c\|_{\max}\le1}}
\left|
\boldsymbol L_{\boldsymbol a}
(\boldsymbol\Sigma-\bar{\boldsymbol\Sigma})
\boldsymbol L_{\boldsymbol c}^\top
\right|
&=
O(M^{-1/2}).
\label{eq:projected-covariance-theorem}
\end{align}
Consequently, for every fixed bounded observation weight $\varphi$,
\begin{align}
\left\|
\operatorname{Cov}(\boldsymbol Y_\varphi\mid\boldsymbol\theta)
-
\boldsymbol\Gamma_\varphi
\right\|_{\rm op}
&=
O(M^{-1/2}).
\label{eq:covariance-theorem}
\end{align}
\end{theorem}

\begin{proof}
Condition on $\boldsymbol\theta$ and suppress this conditioning. Stationarity applied componentwise to
$(\boldsymbol w-\boldsymbol w^*)(\boldsymbol w-\boldsymbol w^*)^\top$ gives
\begin{align}
\boldsymbol J\boldsymbol\Sigma
+
\boldsymbol\Sigma\boldsymbol J^\top
+
\bar{\boldsymbol Q}
&=
-(\boldsymbol E_R+\boldsymbol E_Q).
\label{eq:stationary-second-moment-eq}
\end{align}
Thus,
\begin{align}
\boldsymbol\Sigma-\bar{\boldsymbol\Sigma}
&=
\int_0^\infty
e^{\boldsymbol Jt}
(\boldsymbol E_R+\boldsymbol E_Q)
e^{\boldsymbol J^\top t}dt.
\label{eq:D-integral}
\end{align}
Let
\begin{align}
\delta_R
&=
\sup_{\|\boldsymbol b\|_{\max}\le1}
\frac1M|\boldsymbol b^\top\boldsymbol E_R\boldsymbol b|,
\qquad
\delta_Q
=
\sup_{\|\boldsymbol b\|_{\max}\le1}
\frac1M|\boldsymbol b^\top\boldsymbol E_Q\boldsymbol b|.
\label{eq:delta-RQ}
\end{align}
By A8, $\delta_R=O(M^{-1/2})$ and $\delta_Q=O(M^{-1/2})$. Fix $\boldsymbol E\in\{\boldsymbol E_R,\boldsymbol E_Q\}$ and let
$\delta_E$ denote the corresponding supremum. First suppose that
$\|\boldsymbol a\|_{\max},\|\boldsymbol c\|_{\max}\le1$.
Set $\boldsymbol b_\pm=(\boldsymbol a\pm\boldsymbol c)/2$, so that
$\|\boldsymbol b_\pm\|_{\max}\le1$. Expanding the difference
$\boldsymbol b_+^\top\boldsymbol E\boldsymbol b_+
 -\boldsymbol b_-^\top\boldsymbol E\boldsymbol b_-$ gives
$\tfrac12(\boldsymbol a^\top\boldsymbol E\boldsymbol c
 +\boldsymbol c^\top\boldsymbol E\boldsymbol a)
 =\boldsymbol a^\top\boldsymbol E\boldsymbol c$
by symmetry of $\boldsymbol E$. Each quadratic form in the difference
has absolute value at most $M\delta_E$. Thus
$M^{-1}|\boldsymbol a^\top\boldsymbol E\boldsymbol c|\le2\delta_E$.
For nonzero $\boldsymbol a$ and $\boldsymbol c$, applying this estimate to
$\boldsymbol a/\|\boldsymbol a\|_{\max}$ and
$\boldsymbol c/\|\boldsymbol c\|_{\max}$ gives the following bound,
which also holds if either vector is zero. Therefore, we have
\begin{align}
\frac1M
|\boldsymbol a^\top\boldsymbol E\boldsymbol c|
&\le
2\delta_E
\|\boldsymbol a\|_{\max}\|\boldsymbol c\|_{\max}.
\label{eq:polarization-E}
\end{align}
The dual norm of $\|\cdot\|_{1,M}$ is $M\|\cdot\|_{\max}$, so A6 gives
\begin{align}
\|e^{\boldsymbol J^\top t}\boldsymbol a\|_{\max}
&\le
C_se^{-\gamma_st}\|\boldsymbol a\|_{\max}.
\label{eq:A6-dual}
\end{align}
Combining \eqref{eq:D-integral}--\eqref{eq:A6-dual},
\begin{align}
\sup_{\substack{
\|\boldsymbol a\|_{\max}\le1\\
\|\boldsymbol c\|_{\max}\le1}}
\left|
\boldsymbol L_{\boldsymbol a}\left(\boldsymbol\Sigma-\bar{\boldsymbol\Sigma}\right)\boldsymbol L_{\boldsymbol c}^\top
\right|
&\le
\frac{C_s^2}{\gamma_s}(\delta_R+\delta_Q)
=
O(M^{-1/2}),
\label{eq:projected-D-rate}
\end{align}
which proves \eqref{eq:projected-covariance-theorem}.

For $i=1,\ldots,n-1$, let $\boldsymbol a^{(i)}$ have block
$\boldsymbol a_m^{(i)}=\varphi(\theta_m)\boldsymbol e_i$. Since $n$ is fixed and $\|\boldsymbol a^{(i)}\|_{\max}\le\|\varphi\|_\infty$, \eqref{eq:projected-covariance-theorem} gives
\begin{align}
\left\|
\boldsymbol L_\varphi
(\boldsymbol\Sigma-\bar{\boldsymbol\Sigma})
\boldsymbol L_\varphi^\top
\right\|_{\rm op}
&=
O(M^{-1/2}).
\label{eq:Lphi-D-rate}
\end{align}

By Corollary~\ref{cor:weighted-first-order-moments}, the outer
product of the conditional mean of $\boldsymbol Y_\varphi$
with itself has operator norm $O(M^{-1})$. Since
\begin{align}
\operatorname{Cov}(\boldsymbol Y_\varphi\mid\boldsymbol\theta)
&=
\boldsymbol L_\varphi\boldsymbol\Sigma\boldsymbol L_\varphi^\top
-
\mathbb E[\boldsymbol Y_\varphi\mid\boldsymbol\theta]
\mathbb E[\boldsymbol Y_\varphi\mid\boldsymbol\theta]^\top,
\label{eq:Yphi-cov-identity}
\end{align}
the preceding bound and \eqref{eq:Yphi-cov-identity} prove the covariance
approximation in Theorem~\ref{thm:covariance-approximation}.
\end{proof}

%% file: so-thm-2.tex
We next turn to the refined stationary mean. Define
\begin{align}
\boldsymbol\mu
&=
-\frac12
\boldsymbol J^{-1}\boldsymbol H(\bar{\boldsymbol\Sigma}),
\qquad
\boldsymbol\beta_\varphi
=
\sqrt M\,\boldsymbol L_\varphi\boldsymbol\mu.
\label{eq:mu-beta}
\end{align}

\begin{theorem}[Second-order performance refinement]
\label{thm:performance-refinement}
Assume A1--A9 and fix a bounded observation weight
$\varphi:\Theta\to\mathbb R$. Then
\begin{align}
\left\|
\mathbb E[\boldsymbol X_\varphi(\infty)\mid\boldsymbol\theta]
-\boldsymbol x_\varphi^*
-\frac{\boldsymbol\beta_\varphi}{M}
\right\|
&=o(M^{-1}),
&
\left\|
\operatorname{Cov}(\boldsymbol X_\varphi(\infty)\mid\boldsymbol\theta)
-\frac{\boldsymbol\Gamma_\varphi}{M}
\right\|_{\mathrm{op}}
&=O(M^{-3/2}).
\label{eq:second-order-refinement}
\end{align}
If the sharper $O(M^{-1/2})$ remainder bound in A9 holds, then the
mean error in \eqref{eq:second-order-refinement} is $O(M^{-3/2})$.

Moreover, if $h:\mathbb R^{n-1}\to\mathbb R$ is twice continuously
differentiable on a neighborhood of
${\cal K}_\varphi
=\{\boldsymbol y\in\mathbb R^{n-1}:
\|\boldsymbol y\|_1\leq\|\varphi\|_\infty\}$,
then
\begin{align}
\mathbb E[h(\boldsymbol X_\varphi(\infty))\mid\boldsymbol\theta]
&=
h(\boldsymbol x_\varphi^*)
+\frac{1}{M}
\left(
Dh(\boldsymbol x_\varphi^*)[\boldsymbol\beta_\varphi]
+\frac12
\operatorname{Tr}\!\left(
D^2h(\boldsymbol x_\varphi^*)\boldsymbol\Gamma_\varphi
\right)
\right)
+o(M^{-1}).
\label{eq:second-order-performance}
\end{align}
\end{theorem}

\begin{proof}
Condition on $\boldsymbol\theta$ and suppress this conditioning. Let
\begin{align}
\boldsymbol m
&=\mathbb E[\boldsymbol W(\infty)-\boldsymbol w^*],
&
\boldsymbol D
&=\boldsymbol\Sigma-\overline{\boldsymbol\Sigma}.
\label{eq:thm2-m-D}
\end{align}
By stationarity,
$\boldsymbol J\boldsymbol m+
\mathbb E[\boldsymbol R_F(\boldsymbol W(\infty))]=\boldsymbol 0$.
Using the definition of $\boldsymbol T_F$ in A9 and the linearity of $H$,
\begin{align}
\mathbb E[\boldsymbol R_F(\boldsymbol W(\infty))]
-\frac12H(\overline{\boldsymbol\Sigma})
&=
\mathbb E[\boldsymbol T_F(\boldsymbol W(\infty))]
+\frac12H(\boldsymbol D)
=: \boldsymbol E_H.
\label{eq:thm2-EH}
\end{align}
We first bound $\boldsymbol E_H$. By the Hessian-control part of A9 and
Theorem~\ref{thm:covariance-approximation}, uniformly over $\|\boldsymbol b\|_{\max}\leq1$,
\begin{align}
|\boldsymbol b^\top H(\boldsymbol D)|
&\leq
C_H\sum_{k=1}^{K_H}
\left|
L_{\boldsymbol r_{k,\boldsymbol b}}
\boldsymbol D
L_{\boldsymbol s_{k,\boldsymbol b}}^\top
\right|
=O(M^{-1/2}),
\label{eq:thm2-HD}
\end{align}
where we used
$\|\boldsymbol r_{k,\boldsymbol b}\|_{\max},
\|\boldsymbol s_{k,\boldsymbol b}\|_{\max}\leq C_H$
and rescaled the coefficient vectors when applying Theorem~\ref{thm:covariance-approximation}. Together with
the post-quadratic remainder bound in A9, \eqref{eq:thm2-EH} therefore gives
\begin{align}
\delta_H
&:=
\sup_{\|\boldsymbol b\|_{\max}\leq1}
|\boldsymbol b^\top\boldsymbol E_H|
=o(1).
\label{eq:thm2-deltaH}
\end{align}
If the sharper $O(M^{-1/2})$ remainder bound in A9 holds, then
\begin{align}
\delta_H
&=O(M^{-1/2}).
\label{eq:thm2-deltaH-sharp}
\end{align}

Recall
$\boldsymbol\mu=-\frac12\boldsymbol J^{-1}H(\overline{\boldsymbol\Sigma})$,
so
$\boldsymbol J\boldsymbol\mu+\frac12H(\overline{\boldsymbol\Sigma})
=\boldsymbol0$. Combining this identity with stationarity and
\eqref{eq:thm2-EH} gives
\begin{align}
\boldsymbol m-\boldsymbol\mu
&=
\int_0^\infty e^{\boldsymbol Jt}\boldsymbol E_H\,dt.
\label{eq:thm2-m-minus-mu}
\end{align}
For any $\boldsymbol a\in\mathbb R^{n-1}$ with $\|\boldsymbol a\|=1$,
A6 and duality imply
\begin{align}
\left|
\boldsymbol a^\top L_\varphi
(\boldsymbol m-\boldsymbol\mu)
\right|
&\leq
\int_0^\infty
\left\|
e^{\boldsymbol J^\top t}L_\varphi^\top\boldsymbol a
\right\|_{\max}
\delta_H\,dt
\leq
\frac{C_s\|\varphi\|_\infty}{\gamma_s\sqrt M}\,
\delta_H.
\label{eq:thm2-projected-mean}
\end{align}
Since
$\mathbb E[\boldsymbol X_\varphi(\infty)]-\boldsymbol x_\varphi^*
=M^{-1/2}L_\varphi\boldsymbol m$
and
$\boldsymbol\beta_\varphi/M=M^{-1/2}L_\varphi\boldsymbol\mu$,
\eqref{eq:thm2-deltaH}--\eqref{eq:thm2-projected-mean} yield
\begin{align}
\left\|
\mathbb E[\boldsymbol X_\varphi(\infty)]
-\boldsymbol x_\varphi^*
-\frac{\boldsymbol\beta_\varphi}{M}
\right\|
&=o(M^{-1}).
\label{eq:thm2-mean-result}
\end{align}
Under the sharper version of A9, \eqref{eq:thm2-deltaH-sharp} gives the
sharper bound $O(M^{-3/2})$ in \eqref{eq:thm2-mean-result}.

For the covariance, Theorem~1 gives
\begin{align}
\left\|
L_\varphi(\boldsymbol\Sigma-\overline{\boldsymbol\Sigma})
L_\varphi^\top
\right\|_{\mathrm{op}}
&=O(M^{-1/2}).
\label{eq:thm2-second-moment-cov}
\end{align}
Moreover, by Corollary~\ref{cor:weighted-first-order-moments}, the outer
product of the conditional mean of $\boldsymbol Y_\varphi$
with itself has operator norm $O(M^{-1})$.  Hence
\begin{align}
\left\|
\operatorname{Cov}(\boldsymbol Y_\varphi\mid\boldsymbol\theta)
-\boldsymbol\Gamma_\varphi
\right\|_{\mathrm{op}}
&=O(M^{-1/2}),
&
\left\|
\operatorname{Cov}(\boldsymbol X_\varphi(\infty)\mid\boldsymbol\theta)
-\frac{\boldsymbol\Gamma_\varphi}{M}
\right\|_{\mathrm{op}}
&=O(M^{-3/2}).
\label{eq:thm2-cov-result}
\end{align}

For the performance expansion, Taylor's theorem at $\boldsymbol x_\varphi^*$
gives
\begin{align}
h(\boldsymbol X_\varphi(\infty))
&= h(\boldsymbol x_\varphi^*)
 + Dh(\boldsymbol x_\varphi^*)
   \bigl[\boldsymbol X_\varphi(\infty)-\boldsymbol x_\varphi^*\bigr]
 \notag\\
&\quad + \tfrac12 D^2h(\boldsymbol x_\varphi^*)
 \bigl[\boldsymbol X_\varphi(\infty)-\boldsymbol x_\varphi^*,
       \boldsymbol X_\varphi(\infty)-\boldsymbol x_\varphi^*\bigr]
 + R_h(\boldsymbol X_\varphi(\infty)).
\label{eq:thm2-h-taylor}
\end{align}
By \eqref{eq:thm2-mean-result} and \eqref{eq:thm2-cov-result},
\begin{align}
Dh(\boldsymbol x_\varphi^*)
 \bigl[\mathbb E[\boldsymbol X_\varphi(\infty)]-\boldsymbol x_\varphi^*\bigr]
&= \frac1M Dh(\boldsymbol x_\varphi^*)[\boldsymbol\beta_\varphi]
 + o(M^{-1}), \notag\\
\mathbb E\!\left[D^2h(\boldsymbol x_\varphi^*)
 \bigl[\boldsymbol X_\varphi(\infty)-\boldsymbol x_\varphi^*,
       \boldsymbol X_\varphi(\infty)-\boldsymbol x_\varphi^*\bigr]\right]
&= \frac1M \operatorname{Tr}\!\bigl(
 D^2h(\boldsymbol x_\varphi^*)\boldsymbol\Gamma_\varphi\bigr)
 + o(M^{-1}).
\label{eq:thm2-h-leading}
\end{align}
In the second equality, the outer product of
$\mathbb E[\boldsymbol X_\varphi(\infty)]-\boldsymbol x_\varphi^*$ with
itself has operator norm $O(M^{-2})$ by
Corollary~\ref{cor:weighted-first-order-moments}.

Applying Lemma~\ref{lem:uniform-block-moments} coordinatewise gives
\begin{align}
\mathbb E\bigl\|\boldsymbol X_\varphi(\infty)-\boldsymbol x_\varphi^*\bigr\|^2
&\le \frac{C_\varphi}{M},
&
\mathbb E\bigl\|\boldsymbol X_\varphi(\infty)-\boldsymbol x_\varphi^*\bigr\|^4
&\le \frac{C_\varphi}{M^2}.
\label{eq:thm2-zeta-moments}
\end{align}
Let
\begin{align}
\omega_h(r)
&= \sup_{\substack{\boldsymbol y,\boldsymbol z\in{\cal K}_\varphi\\
                    \|\boldsymbol y-\boldsymbol z\|\le r}}
   \bigl\|D^2h(\boldsymbol y)-D^2h(\boldsymbol z)\bigr\|_{\mathrm{op}},
&
B_h
&= \sup_{\boldsymbol y\in{\cal K}_\varphi}
   \|D^2h(\boldsymbol y)\|_{\mathrm{op}}.
\label{eq:thm2-h-modulus}
\end{align}
The Taylor remainder satisfies
$|R_h(\boldsymbol X_\varphi(\infty))|
 \le \tfrac12\omega_h\bigl(\|\boldsymbol X_\varphi(\infty)-\boldsymbol x_\varphi^*\|\bigr)
      \|\boldsymbol X_\varphi(\infty)-\boldsymbol x_\varphi^*\|^2$.
For fixed $\eta>0$, split the expectation over
$\{\|\boldsymbol X_\varphi(\infty)-\boldsymbol x_\varphi^*\|\le\eta\}$
and its complement. Using $\omega_h(r)\le 2B_h$ and
\eqref{eq:thm2-zeta-moments} gives
\begin{align}
M\mathbb E\bigl|R_h(\boldsymbol X_\varphi(\infty))\bigr|
&\le \frac{C_\varphi}{2}\omega_h(\eta)
 + \frac{C_\varphi B_h}{M\eta^2}.
\label{eq:thm2-h-remainder}
\end{align}
Since $D^2h$ is uniformly continuous on ${\cal K}_\varphi$,
$\omega_h(\eta)\to0$ as $\eta\downarrow0$. Taking $M\to\infty$ and then
$\eta\downarrow0$ in \eqref{eq:thm2-h-remainder} yields
\begin{align}
\mathbb E\bigl|R_h(\boldsymbol X_\varphi(\infty))\bigr|
&= o(M^{-1}).
\label{eq:thm2-h-remainder-small}
\end{align}
Taking expectations in \eqref{eq:thm2-h-taylor} and using
\eqref{eq:thm2-h-leading} and \eqref{eq:thm2-h-remainder-small} proves the
performance expansion.
\end{proof}

%% file: po2.tex
\section{Application: Heterogeneous Power-of-Two Choices with $B=2$}
\label{sec:po2-b2}
We consider a heterogeneous power-of-two-choices system with $M$ servers and buffer size $B=2$, including the job in service. Server $m$ has service rate $\mu_m$, where
\begin{align}
0<\underline\mu\leq \mu_m\leq\overline\mu<\infty,
\qquad m=1,\ldots,M.
\label{eq:po2-mu-bound}
\end{align}
The arrival-rate parameter $\lambda>0$ is fixed independently of $M$. Jobs arrive according to a Poisson process of rate $\lambda M$. Each arrival samples two servers independently and uniformly with replacement and is routed to the server with the shorter queue; ties are broken uniformly. An arrival is discarded if both sampled servers contain two jobs. Throughout this section, we condition on the realized service-rate vector $\boldsymbol\mu=(\mu_1,\ldots,\mu_M)^\top$.

Let $Q_m(t)\in\{0,1,2\}$ be the queue length of server $m$ and use the following state representation
\begin{align}
\boldsymbol W_m(t)
&=
\begin{pmatrix}
W_{m,1}(t)\\ W_{m,2}(t)
\end{pmatrix}
=
\begin{pmatrix}
\mathbf{1}\{Q_m(t)\geq1\}\\ \mathbf{1}\{Q_m(t)=2\}
\end{pmatrix}
\in
{\cal S}
=
\left\{
\begin{pmatrix}0\\0\end{pmatrix},
\begin{pmatrix}1\\0\end{pmatrix},
\begin{pmatrix}1\\1\end{pmatrix}
\right\}.
\label{eq:po2-tail-state}
\end{align}
This is an invertible linear transformation of the reduced one-hot representation. Specifically, if
$\widetilde{\boldsymbol W}_m=(\mathbf{1}\{Q_m=1\},\mathbf{1}\{Q_m=2\})^\top$, then
\begin{align}
\boldsymbol W_m
&=
\begin{pmatrix}1&1\\0&1\end{pmatrix}
\widetilde{\boldsymbol W}_m.
\label{eq:po2-tail-transform}
\end{align}
Since the transformation and its inverse have fixed operator norms, the
assumptions and conclusions of Sections~\ref{sec:first-order}
and~\ref{sec:second-order} carry over with constants independent of $M$.
In tail coordinates, the single-server mean-field state lies in the triangle
in~\eqref{eq:po2-domain}, the empirical tail occupancy lies in
the triangle in~\eqref{eq:po2-K}, and the weighted tail
occupancy lies in the box in~\eqref{eq:po2-Kphi}.

Define the empirical tail occupancies
\begin{align}
X_i(\boldsymbol W)
&=\frac1M\sum_{m=1}^M W_{m,i},
\qquad i=1,2.
\label{eq:po2-tail-occupancy}
\end{align}
This model fits the framework of Section~\ref{sec:models} with $\theta_m=\mu_m$, $\chi\equiv1$, and $\Psi$ equal to the identity map on $\mathbb R^2$, so the interaction field is
$\boldsymbol Z(\boldsymbol W)=(X_1(\boldsymbol W),X_2(\boldsymbol W))^\top$. The transition rates of server $m$ are
\begin{align}
\begin{pmatrix}0\\0\end{pmatrix}
&\longrightarrow
\begin{pmatrix}1\\0\end{pmatrix}
&&\text{at rate }\lambda(1+X_1),
&
\begin{pmatrix}1\\0\end{pmatrix}
&\longrightarrow
\begin{pmatrix}1\\1\end{pmatrix}
&&\text{at rate }\lambda(X_1+X_2),
\nonumber\\
\begin{pmatrix}1\\0\end{pmatrix}
&\longrightarrow
\begin{pmatrix}0\\0\end{pmatrix}
&&\text{at rate }\mu_m,
&
\begin{pmatrix}1\\1\end{pmatrix}
&\longrightarrow
\begin{pmatrix}1\\0\end{pmatrix}
&&\text{at rate }\mu_m.
\label{eq:po2-transition-rates}
\end{align}
Because $\lambda>0$, $\mu_m>0$, and sampling is with replacement, the finite CTMC is irreducible and therefore has a unique stationary distribution. For
$\boldsymbol w_m=(w_{m,1},w_{m,2})^\top$ and
$x_i(\boldsymbol w)=M^{-1}\sum_m w_{m,i}$, the mean-field drift is
\begin{align}
F_{m,1}(\boldsymbol w)
&=\lambda(1+x_1)(1-w_{m,1})-\mu_m(w_{m,1}-w_{m,2}),
\nonumber\\
F_{m,2}(\boldsymbol w)
&=\lambda(x_1+x_2)(w_{m,1}-w_{m,2})-\mu_m w_{m,2}.
\label{eq:po2-drift}
\end{align}
The admissible state space of one mean-field block is
\begin{align}
{\cal D}_{\rm po2}
&=\left\{(w_1,w_2)^\top\in\mathbb R^2:0\leq w_2\leq w_1\leq1\right\}.
\label{eq:po2-domain}
\end{align}
For this application, we take the full admissible mean-field domain in Section~\ref{sec:models} to be
${\cal D}={\cal D}_{\rm po2}^{\times M}$, which is compact, convex, and forward invariant under \eqref{eq:po2-drift}.

For fixed $(x_1,x_2)$, let
\begin{align}
a&=\lambda(1+x_1),
&
b&=\lambda(x_1+x_2),
&
D_m&=\mu_m^2+\mu_m a+ab.
\label{eq:po2-abD}
\end{align}
The stationary tail probabilities of the corresponding three-state birth--death chain are
\begin{align}
w_{m,1}
&=\frac{a(\mu_m+b)}{D_m},
&
w_{m,2}
&=\frac{ab}{D_m}.
\label{eq:po2-local-fixed-point}
\end{align}
Consequently, \eqref{eq:po2-local-fixed-point} defines a continuous map from
\begin{align}
{\cal K}
&=\left\{(x_1,x_2)^\top\in\mathbb R^2:0\leq x_2\leq x_1\leq1\right\}
\label{eq:po2-K}
\end{align}
into itself by averaging over $m$. Brouwer's fixed-point theorem therefore gives at least one equilibrium $\boldsymbol w^*$. Its uniqueness follows from the incremental stability proved below.

\subsection{Verification of S1--S2}
\label{subsec:po2s1-s2}
We first verify S1. For $\boldsymbol v=(v_1,v_2)^\top$ and
$\boldsymbol x=(x_1,x_2)^\top$, define
\begin{align}
\boldsymbol A(\boldsymbol x;\mu)
&=
\begin{pmatrix}
-\lambda(1+x_1)-\mu & \mu\\
\lambda(x_1+x_2) & -\lambda(x_1+x_2)-\mu
\end{pmatrix},
&
\boldsymbol a(\boldsymbol x;\mu)
&=
\begin{pmatrix}
\lambda(1+x_1)\\0
\end{pmatrix}.
\label{eq:po2-Aa}
\end{align}
Then the drift in \eqref{eq:po2-drift} can be written as
\begin{align}
f(\boldsymbol v,\boldsymbol x;\mu)
&=\boldsymbol A(\boldsymbol x;\mu)\boldsymbol v+\boldsymbol a(\boldsymbol x;\mu).
\label{eq:po2-local-affine}
\end{align}
Thus the new local approximately affine condition S1 holds with residual
$\boldsymbol r\equiv\boldsymbol0$. In fact, \eqref{eq:po2-local-affine} holds globally on
${\cal D}_{\rm po2}\times{\cal K}$, so any fixed local radius $\eta_A>0$ allowed by the interaction domain may be used. Moreover, $\boldsymbol A$ and $\boldsymbol a$ are affine in $\boldsymbol x$, their first derivatives are uniformly bounded under \eqref{eq:po2-mu-bound}, and their second derivatives vanish identically. Hence one may take $\omega_A\equiv\omega_r\equiv0$, so both S1 and its sharper version hold. In particular, the drift is exactly quadratic in the full mean-field state and the post-quadratic remainder $\boldsymbol T_F$ vanishes identically.

We next verify S2. Each transition changes exactly one tail coordinate by $\pm1$, so
\begin{align}
\boldsymbol Q(\boldsymbol w)
&=\operatorname{diag}\bigl(q_{1,1}(\boldsymbol w),q_{1,2}(\boldsymbol w),\ldots,q_{M,1}(\boldsymbol w),q_{M,2}(\boldsymbol w)\bigr),
\label{eq:po2-Q}
\end{align}
where
\begin{align}
q_{m,1}(\boldsymbol w)
&=\lambda(1-w_{m,1})(1+x_1)+\mu_m(w_{m,1}-w_{m,2}),
\nonumber\\
q_{m,2}(\boldsymbol w)
&=\lambda(w_{m,1}-w_{m,2})(x_1+x_2)+\mu_m w_{m,2}.
\label{eq:po2-q12}
\end{align}
Fix $\boldsymbol b$ with $\|\boldsymbol b\|_{\max}\leq1$, and define
$\boldsymbol e^{(1)},\boldsymbol e^{(+)}\in\mathbb R^{2M}$ by
\begin{align}
\boldsymbol e_m^{(1)}
&=(1,0)^\top,
&
\boldsymbol e_m^{(+)}
&=(1,1)^\top.
\label{eq:po2-e-vectors}
\end{align}
Also define
\begin{align}
\boldsymbol a_{\boldsymbol b,m}
&=
\begin{pmatrix}
 b_{m,1}^2\bigl(\mu_m-\lambda(1+x_1^*)\bigr)
 +b_{m,2}^2\lambda(x_1^*+x_2^*)\\
 -b_{m,1}^2\mu_m
 +b_{m,2}^2\bigl(\mu_m-\lambda(x_1^*+x_2^*)\bigr)
\end{pmatrix},
\nonumber\\
\boldsymbol d_{\boldsymbol b,m}^{(1)}
&=(b_{m,1}^2,0)^\top,
&
\boldsymbol d_{\boldsymbol b,m}^{(2)}
&=(b_{m,2}^2,-b_{m,2}^2)^\top,
\nonumber\\
\rho_{1,\boldsymbol b}
&=\frac{\lambda}{M}\sum_{m=1}^M b_{m,1}^2(1-w_{m,1}^*),
&
\rho_{2,\boldsymbol b}
&=\frac{\lambda}{M}\sum_{m=1}^M b_{m,2}^2(w_{m,1}^*-w_{m,2}^*).
\label{eq:po2-S2-vectors}
\end{align}
Expanding \eqref{eq:po2-q12} around $\boldsymbol w^*$ gives the exact identity
\begin{align}
V_{\boldsymbol b}(\boldsymbol w)
&=V_{\boldsymbol b}(\boldsymbol w^*)
+\bar{\boldsymbol x}_{\boldsymbol a_{\boldsymbol b}}(\boldsymbol w)
+\rho_{1,\boldsymbol b}\bar{\boldsymbol x}_{\boldsymbol e^{(1)}}(\boldsymbol w)
+\rho_{2,\boldsymbol b}\bar{\boldsymbol x}_{\boldsymbol e^{(+)}}(\boldsymbol w)
\nonumber\\
&\quad
-\lambda \bar{\boldsymbol x}_{\boldsymbol d_{\boldsymbol b}^{(1)}}(\boldsymbol w)
 \bar{\boldsymbol x}_{\boldsymbol e^{(1)}}(\boldsymbol w)
+\lambda \bar{\boldsymbol x}_{\boldsymbol d_{\boldsymbol b}^{(2)}}(\boldsymbol w)
 \bar{\boldsymbol x}_{\boldsymbol e^{(+)}}(\boldsymbol w).
\label{eq:po2-S2-identity}
\end{align}
Thus S2 has a direct finite weighted-occupancy representation. Specifically, take $K_Q=5$ with
\begin{align}
\boldsymbol c_{1,\boldsymbol b}
&=\boldsymbol a_{\boldsymbol b},
&
\boldsymbol c_{2,\boldsymbol b}
&=\boldsymbol e^{(1)},
&
\boldsymbol c_{3,\boldsymbol b}
&=\boldsymbol e^{(+)},
&
\boldsymbol c_{4,\boldsymbol b}
&=\boldsymbol d_{\boldsymbol b}^{(1)},
&
\boldsymbol c_{5,\boldsymbol b}
&=\boldsymbol d_{\boldsymbol b}^{(2)}.
\label{eq:po2-S2-c}
\end{align}
For $\boldsymbol z=(z_1,\ldots,z_5)^\top$, define
\begin{align}
\widehat V_{\boldsymbol b}(\boldsymbol z)
&=V_{\boldsymbol b}(\boldsymbol w^*)+z_1
+\rho_{1,\boldsymbol b}z_2+\rho_{2,\boldsymbol b}z_3
-\lambda z_4z_2+\lambda z_5z_3.
\label{eq:po2-Vhat}
\end{align}
Then \eqref{eq:po2-S2-identity} is exactly
\begin{align}
V_{\boldsymbol b}(\boldsymbol w)
&=\widehat V_{\boldsymbol b}\!\left(
\bar{\boldsymbol x}_{\boldsymbol c_{1,\boldsymbol b}}(\boldsymbol w),\ldots,
\bar{\boldsymbol x}_{\boldsymbol c_{5,\boldsymbol b}}(\boldsymbol w)
\right).
\label{eq:po2-S2-representation}
\end{align}
All vectors in \eqref{eq:po2-S2-c} have uniformly bounded block-max norms, and
$|\rho_{1,\boldsymbol b}|+|\rho_{2,\boldsymbol b}|\leq2\lambda$. Moreover, for every
$\boldsymbol w\in{\cal D}_{\rm po2}^{\times M}$,
\begin{align}
q_{m,1}(\boldsymbol w)+q_{m,2}(\boldsymbol w)
&\leq
2\lambda\bigl((1-w_{m,1})+(w_{m,1}-w_{m,2})\bigr)+\mu_m w_{m,1}
\leq2\lambda+\overline\mu.
\label{eq:po2-q-bound}
\end{align}
Therefore, since $b_{m,i}^2\leq1$,
\begin{align}
0\leq V_{\boldsymbol b}(\boldsymbol w^*)
&\leq\frac1M\sum_{m=1}^M\bigl(q_{m,1}(\boldsymbol w^*)+q_{m,2}(\boldsymbol w^*)\bigr)
\leq2\lambda+\overline\mu.
\label{eq:po2-Vstar-bound}
\end{align}
Finally, $\widehat V_{\boldsymbol b}$ in \eqref{eq:po2-Vhat} is a quadratic polynomial. Hence, for any fixed $\eta_Q>0$, $\widehat V_{\boldsymbol b}$ and its first two derivatives are uniformly bounded on
$\{\|\boldsymbol z\|_\infty\leq\eta_Q\}$, with constants depending only on
$\lambda$, $\overline\mu$, and $\eta_Q$. Therefore S2 holds uniformly in
$M$, $\boldsymbol\mu$, and $\boldsymbol b$.

\subsection{Verification of A1--A6}

Set
\begin{align}
\vartheta
&=\min\left\{1,\frac{\underline\mu}{4\lambda}\right\},
&
c&=1+\vartheta,
&
\beta&=\frac{\vartheta\underline\mu}{1+\vartheta}>0.
\label{eq:po2-vartheta}
\end{align}

\begin{lemma}
Under \eqref{eq:po2-mu-bound}, the heterogeneous power-of-two-choices model with $B=2$ satisfies A1--A6 for every fixed $\lambda>0$.
\label{lem:po2-A1-A6}
\end{lemma}

\begin{proof}
A1 follows because the total accepted-arrival rate is at most $\lambda M$ and the total service rate is at most $\overline\mu M$. A2 holds because every transition changes the state of one server. A3 follows from \eqref{eq:po2-drift}, whose entries are quadratic polynomials with derivatives uniformly bounded on ${\cal D}_{\rm po2}^{\times M}$.

We first establish a weighted $1$-norm contraction. Let
$\boldsymbol J(\boldsymbol w)=D\boldsymbol F(\boldsymbol w)$ and let
$\boldsymbol p\in\mathbb R^{2M}$ have every block equal to $(1,c)^\top$. Since
$\boldsymbol J(\boldsymbol w)$ is Metzler, it suffices to bound its weighted column sums. Summing the two drift coordinates with weights $1$ and $c$ gives
\begin{align}
\sum_{m=1}^M\bigl(F_{m,1}(\boldsymbol w)+cF_{m,2}(\boldsymbol w)\bigr)
&=\lambda M\bigl(1+\vartheta x_1^2-cx_2^2\bigr)
-\sum_{m=1}^M\mu_m w_{m,1}
-\vartheta\sum_{m=1}^M\mu_m w_{m,2}.
\label{eq:po2-weighted-sum}
\end{align}
Differentiating \eqref{eq:po2-weighted-sum} with respect to $w_{j,1}$ and $w_{j,2}$ yields
\begin{align}
(\boldsymbol p^\top\boldsymbol J(\boldsymbol w))_{j,1}
&=-\mu_j+2\lambda\vartheta x_1
\leq-\underline\mu+2\lambda\vartheta
\leq-\frac{\underline\mu}{2}
\leq-\beta,
\nonumber\\
(\boldsymbol p^\top\boldsymbol J(\boldsymbol w))_{j,2}
&=-\vartheta\mu_j-2c\lambda x_2
\leq-\vartheta\underline\mu
=-\beta c.
\label{eq:po2-column-sum}
\end{align}
For any two mean-field trajectories $\boldsymbol w(t)$ and
$\widetilde{\boldsymbol w}(t)$, their difference
$\boldsymbol u(t)=\boldsymbol w(t)-\widetilde{\boldsymbol w}(t)$ satisfies
\begin{align}
\dot{\boldsymbol u}(t)
&=\boldsymbol J(t)\boldsymbol u(t),
&
\boldsymbol J(t)
&=\int_0^1\boldsymbol J\!\left(
\widetilde{\boldsymbol w}(t)+s\boldsymbol u(t)
\right)ds.
\label{eq:po2-difference-J}
\end{align}
The matrix $\boldsymbol J(t)$ is Metzler and satisfies \eqref{eq:po2-column-sum}. Therefore, with
\begin{align}
V(t)
&=\frac1M\sum_{m=1}^M\bigl(|u_{m,1}(t)|+c|u_{m,2}(t)|\bigr),
\label{eq:po2-V}
\end{align}
the upper-right Dini derivative satisfies
\begin{align}
D^+V(t)
&\leq-\beta V(t),
&
V(t)&\leq e^{-\beta t}V(0).
\label{eq:po2-V-contraction}
\end{align}
We next prove A4 in the normalized Euclidean norm. Set
$\overline u_i=M^{-1}\sum_m u_{m,i}$. Direct subtraction of the two mean-field systems gives
\begin{align}
\dot u_{m,1}
&=-\bigl(\mu_m+\lambda(1+x_1)\bigr)u_{m,1}+\mu_m u_{m,2}
+\lambda(1-\widetilde w_{m,1})\overline u_1,
\nonumber\\
\dot u_{m,2}
&=\lambda(x_1+x_2)u_{m,1}
-\bigl(\mu_m+\lambda(x_1+x_2)\bigr)u_{m,2}
+\lambda(\widetilde w_{m,1}-\widetilde w_{m,2})(\overline u_1+\overline u_2),
\label{eq:po2-difference-block}
\end{align}
where $x_i=M^{-1}\sum_m w_{m,i}$ refers to the first trajectory. The local homogeneous matrix in \eqref{eq:po2-difference-block} is
\begin{align}
\boldsymbol B_m(t)
&=
\begin{pmatrix}
-\mu_m-\lambda(1+x_1(t)) & \mu_m\\
\lambda(x_1(t)+x_2(t)) & -\mu_m-\lambda(x_1(t)+x_2(t))
\end{pmatrix}.
\label{eq:po2-Bm}
\end{align}
Its row sums are $-\lambda(1+x_1(t))\leq-\lambda$ and $-\mu_m\leq-\underline\mu$. Hence, with $a_0=\min\{\lambda,\underline\mu\}$, the evolution operator $\boldsymbol\Phi_m(t,s)$ of
$\dot{\boldsymbol v}=\boldsymbol B_m(t)\boldsymbol v$ satisfies
\begin{align}
\|\boldsymbol\Phi_m(t,s)\|_\infty
&\leq e^{-a_0(t-s)},
&
\|\boldsymbol\Phi_m(t,s)\|_{\mathrm{op}}
&\leq\sqrt2e^{-a_0(t-s)}.
\label{eq:po2-Phi}
\end{align}
The forcing term in \eqref{eq:po2-difference-block} has Euclidean norm at most $2\lambda V(t)$. Therefore, by variation of constants, Minkowski's inequality, \eqref{eq:po2-V-contraction}, and
$V(0)\leq\sqrt{1+c^2}\|\boldsymbol u(0)\|_M$,
\begin{align}
\|\boldsymbol u(t)\|_M
&\leq
\sqrt2e^{-a_0t}\|\boldsymbol u(0)\|_M
+2\sqrt2\lambda\sqrt{1+c^2}\|\boldsymbol u(0)\|_M
\int_0^t e^{-a_0(t-s)}e^{-\beta s}ds
\nonumber\\
&\leq C_0e^{-\gamma_0t}\|\boldsymbol u(0)\|_M,
\qquad
\gamma_0=\frac12\min\{a_0,\beta\},
\label{eq:po2-A4}
\end{align}
where $C_0<\infty$ depends only on $\lambda$, $\underline\mu$, and $\overline\mu$. Taking
$\widetilde{\boldsymbol w}(t)=\boldsymbol w^*$ proves A4 and also proves uniqueness of the equilibrium.

For A5, fix a feasible CTMC transition
$\boldsymbol W\to\boldsymbol W+\boldsymbol\Delta_r$ and a point
$\widehat{\boldsymbol w}=\boldsymbol W+s\boldsymbol\Delta_r$, $s\in[0,1]$, on the corresponding jump segment. Define the first and second variations
\begin{align}
\boldsymbol v(t)
&=D_{\widehat{\boldsymbol w}}\boldsymbol w(t;\widehat{\boldsymbol w})[\boldsymbol\Delta_r],
&
\boldsymbol r(t)
&=D^2_{\widehat{\boldsymbol w}}\boldsymbol w(t;\widehat{\boldsymbol w})[\boldsymbol\Delta_r,\boldsymbol\Delta_r],
\label{eq:po2-variations}
\end{align}
and
\begin{align}
V_1(t)
&=\frac1M\sum_{m=1}^M\bigl(|v_{m,1}(t)|+c|v_{m,2}(t)|\bigr),
&
V_2(t)
&=\frac1M\sum_{m=1}^M\bigl(|r_{m,1}(t)|+c|r_{m,2}(t)|\bigr).
\label{eq:po2-V12}
\end{align}
The first variation satisfies the linearized equation along the trajectory, so \eqref{eq:po2-column-sum} gives
\begin{align}
V_1(t)
&\leq e^{-\beta t}V_1(0)
\leq\frac{c}{M}e^{-\beta t}.
\label{eq:po2-V1}
\end{align}
The Hessian forcing in the second-variation equation is
\begin{align}
D^2F_{m,1}(\boldsymbol w)[\boldsymbol v,\boldsymbol v]
&=-2\lambda v_{m,1}\overline v_1,
&
D^2F_{m,2}(\boldsymbol w)[\boldsymbol v,\boldsymbol v]
&=2\lambda(\overline v_1+\overline v_2)(v_{m,1}-v_{m,2}),
\label{eq:po2-Hessian-force}
\end{align}
where $\overline v_i=M^{-1}\sum_m v_{m,i}$. Consequently,
\begin{align}
\frac1M\sum_{m=1}^M
\left(
|D^2F_{m,1}[\boldsymbol v,\boldsymbol v]|
+c|D^2F_{m,2}[\boldsymbol v,\boldsymbol v]|
\right)
&\leq2\lambda(1+c)V_1(t)^2.
\label{eq:po2-Hessian-force-bound}
\end{align}
Since $\boldsymbol r(0)=\boldsymbol0$, \eqref{eq:po2-V1}--\eqref{eq:po2-Hessian-force-bound} imply
\begin{align}
V_2(t)
&\leq2\lambda(1+c)\int_0^t e^{-\beta(t-s)}V_1(s)^2ds
\leq\frac{2\lambda(1+c)c^2}{\beta M^2}e^{-\beta t}.
\label{eq:po2-V2}
\end{align}
Because $\|\boldsymbol v(t)\|_{1,M}\leq V_1(t)$ and
$\|\boldsymbol r(t)\|_{1,M}\leq V_2(t)$,
\begin{align}
\int_0^\infty\|\boldsymbol v(t)\|_{1,M}^2dt
&\leq\frac{c^2}{2\beta M^2},
&
\int_0^\infty\|\boldsymbol r(t)\|_{1,M}dt
&\leq\frac{2\lambda(1+c)c^2}{\beta^2M^2},
\label{eq:po2-A5}
\end{align}
which proves A5 uniformly along every feasible jump segment.

Finally, the equilibrium Jacobian $\boldsymbol J=D\boldsymbol F(\boldsymbol w^*)$ satisfies the same weighted column-sum bound. Define
\begin{align}
\|\boldsymbol v\|_{c,1,M}
&=\frac1M\sum_{m=1}^M\bigl(|v_{m,1}|+c|v_{m,2}|\bigr).
\label{eq:po2-weighted-norm}
\end{align}
Then
\begin{align}
\|\boldsymbol v\|_{1,M}
&\leq\|\boldsymbol v\|_{c,1,M}
\leq\sqrt{1+c^2}\|\boldsymbol v\|_{1,M}.
\label{eq:po2-norm-equivalence}
\end{align}
Hence, using the weighted contraction for the linearized system,
\begin{align}
\|e^{\boldsymbol Jt}\boldsymbol v\|_{1,M}
&\leq\|e^{\boldsymbol Jt}\boldsymbol v\|_{c,1,M}
\leq e^{-\beta t}\|\boldsymbol v\|_{c,1,M}
\leq\sqrt{1+c^2}e^{-\beta t}\|\boldsymbol v\|_{1,M},
\qquad t\geq0,
\label{eq:po2-A6}
\end{align}
so A6 holds.
\end{proof}

\subsection{Verification of A7}

Define the single-server level function
\begin{align}
\ell(0)&=0,
&
\ell(1)&=1,
&
\ell(2)&=2+\vartheta,
\label{eq:po2-level-function}
\end{align}
and the additive metric on $\{0,1,2\}^M$ by
\begin{align}
d(\boldsymbol q,\widetilde{\boldsymbol q})
&=\sum_{m=1}^M|\ell(q_m)-\ell(\widetilde q_m)|.
\label{eq:po2-metric}
\end{align}
Two states are neighboring in the sense of A7 precisely when they differ by one queue level at one server. Such an edge has length $1$ for a $0\leftrightarrow1$ move and $c=1+\vartheta$ for a $1\leftrightarrow2$ move. Moreover, any two states can be connected by changing the queue lengths one server at a time along $0\leftrightarrow1\leftrightarrow2$. Therefore, \eqref{eq:po2-metric} is the corresponding shortest-path metric and
\begin{align}
1\leq d(\boldsymbol q,\widetilde{\boldsymbol q})
&\leq c,
\qquad \boldsymbol q\sim\widetilde{\boldsymbol q},
&
\sup_{\boldsymbol q,\widetilde{\boldsymbol q}\in\{0,1,2\}^M}
d(\boldsymbol q,\widetilde{\boldsymbol q})
&\leq(2+\vartheta)M.
\label{eq:po2-metric-bounds}
\end{align}

\begin{lemma}[Local coupling contraction]
Under \eqref{eq:po2-mu-bound}, A7 holds with
$D_e=c$ and $\kappa_c=\beta$.
\label{lem:po2-A7}
\end{lemma}

\begin{proof}
By \eqref{eq:po2-metric-bounds}, it remains only to verify the local coupled-generator inequality in A7. By symmetry, there are two types of neighboring pairs. Couple arrivals by using the same arrival clock, the same ordered sampled pair, and the same tie-breaking random variable in both systems. Couple services by using the same rate-$\mu_j$ potential-service clock at each server $j$ in both systems.

Suppose first that the two states differ only by $q_m=0$ and $\widetilde q_m=1$. A service completion at server $m$ coalesces the two states and decreases the distance by $1$. If an arrival does not sample server $m$, the two systems make the same routing decision and the distance is unchanged. If an arrival samples server $m$, the post-arrival discrepancy is either another $0\leftrightarrow1$ discrepancy or a $1\leftrightarrow2$ discrepancy, so the distance can increase by at most $\vartheta$. The rate at which an arrival samples server $m$ at least once is
\begin{align}
\lambda M\left[1-\left(1-\frac1M\right)^2\right]
&=\lambda\left(2-\frac1M\right)
\leq2\lambda.
\label{eq:po2-sample-rate}
\end{align}
Hence the coupled generator satisfies
\begin{align}
({\cal G}^c d)(\boldsymbol q,\widetilde{\boldsymbol q})
&\leq-\mu_m+2\lambda\vartheta
\leq-\frac{\underline\mu}{2}
\leq-\beta d(\boldsymbol q,\widetilde{\boldsymbol q}).
\label{eq:po2-coupling-01}
\end{align}
Suppose next that the states differ only by $q_m=1$ and $\widetilde q_m=2$. Arrivals cannot increase the distance under the synchronous routing coupling: they either leave the discrepancy unchanged, move the same $1\leftrightarrow2$ discrepancy to another sampled server, or coalesce the two systems. A service completion at server $m$ changes the discrepancy from $1\leftrightarrow2$ to $0\leftrightarrow1$, decreasing the distance by $\vartheta$. Therefore,
\begin{align}
({\cal G}^c d)(\boldsymbol q,\widetilde{\boldsymbol q})
&\leq-\mu_m\vartheta
\leq-\frac{\vartheta\underline\mu}{1+\vartheta}
d(\boldsymbol q,\widetilde{\boldsymbol q})
=-\beta d(\boldsymbol q,\widetilde{\boldsymbol q}).
\label{eq:po2-coupling-12}
\end{align}
Thus the local contraction in A7 holds for every neighboring pair. Together with \eqref{eq:po2-metric-bounds} and the shortest-path property of \eqref{eq:po2-metric}, this verifies A7.
\end{proof}

\subsection{Second-Order Consequences}

For a fixed bounded weight $\varphi$, define
\begin{align}
\boldsymbol Y_\varphi
&=\sqrt M\left(\boldsymbol X_\varphi(\boldsymbol W(\infty))-\boldsymbol x_\varphi^*\right),
&
\boldsymbol L_\varphi
&=\frac1{\sqrt M}
\left(\varphi(\mu_1)\boldsymbol I_2\ \cdots\ \varphi(\mu_M)\boldsymbol I_2\right).
\label{eq:po2-YL}
\end{align}
Let $\overline{\boldsymbol\Sigma}$ solve
\begin{align}
\boldsymbol J\overline{\boldsymbol\Sigma}
+\overline{\boldsymbol\Sigma}\boldsymbol J^\top
+\boldsymbol Q(\boldsymbol w^*)
&=\boldsymbol0.
\label{eq:po2-Lyapunov}
\end{align}

\begin{corollary}[Second-order approximation for heterogeneous power-of-two choices with $B=2$]
Under \eqref{eq:po2-mu-bound}, for every fixed $\lambda>0$, conditions A1--A9 hold. For every fixed bounded observation weight $\varphi$, define
\begin{align}
\boldsymbol\Gamma_\varphi
&=\boldsymbol L_\varphi\overline{\boldsymbol\Sigma}\boldsymbol L_\varphi^\top,
&
\boldsymbol\beta_\varphi
&=-\frac{\sqrt M}{2}\boldsymbol L_\varphi\boldsymbol J^{-1}H(\overline{\boldsymbol\Sigma}).
\label{eq:po2-Gamma-beta}
\end{align}
Then
\begin{align}
\left\|\operatorname{Cov}(\boldsymbol Y_\varphi\mid\boldsymbol\mu)
-\boldsymbol\Gamma_\varphi\right\|_{\mathrm{op}}
&=O(M^{-1/2}),
&
\left\|
\mathbb E[\boldsymbol X_\varphi(\boldsymbol W(\infty))\mid\boldsymbol\mu]
-\boldsymbol x_\varphi^*
-\frac{\boldsymbol\beta_\varphi}{M}
\right\|
&=O(M^{-3/2}).
\label{eq:po2-second-order}
\end{align}
Moreover, let
\begin{align}
{\cal K}_\varphi
&=\left\{\boldsymbol y\in\mathbb R^2:|y_i|\leq\|\varphi\|_\infty,\ i=1,2\right\}.
\label{eq:po2-Kphi}
\end{align}
For every $h:\mathbb R^2\to\mathbb R$ that is twice continuously differentiable on a neighborhood of ${\cal K}_\varphi$,
\begin{align}
\mathbb E[h(\boldsymbol X_\varphi(\boldsymbol W(\infty)))\mid\boldsymbol\mu]
&=h(\boldsymbol x_\varphi^*)
+\frac1M\left(
Dh(\boldsymbol x_\varphi^*)[\boldsymbol\beta_\varphi]
+\frac12\operatorname{Tr}\!\left(
D^2h(\boldsymbol x_\varphi^*)\boldsymbol\Gamma_\varphi
\right)
\right)
+o(M^{-1}).
\label{eq:po2-performance-expansion}
\end{align}
\end{corollary}

\begin{proof}
Lemmas \ref{lem:po2-A1-A6} and \ref{lem:po2-A7} establish A1-A7. Section~\ref{subsec:po2s1-s2} verifies S1--S2, so Lemma \ref{lem:structural-A8-A9} gives A8--A9. Hence the covariance and performance conclusions follow from Theorems~1--2. Moreover, here $\boldsymbol r\equiv\boldsymbol0$ and $D^2\boldsymbol A=D^2\boldsymbol a=0$, so the post-quadratic drift remainder vanishes identically and the sharper A9 bound holds. The $O(M^{-3/2})$ refined-mean rate therefore follows as well.
\end{proof}